\documentclass[10pt,a4paper]{article}
\usepackage[english]{babel}
\usepackage[fleqn]{amsmath}
\usepackage{courier}
\usepackage{amsthm}
\usepackage{amsmath}
\usepackage[letterspace=150]{microtype}
\usepackage{lipsum}
\usepackage{amsfonts}
\usepackage{amssymb}
\allowdisplaybreaks
\usepackage{amsfonts}
\usepackage{amssymb}
\usepackage{graphicx}
\usepackage{xcolor}
\usepackage{enumitem}
\usepackage{array}
\usepackage{multirow}
\usepackage{tablists}
\usepackage{cancel}
\numberwithin{equation}{section}

\newtheorem{Theorem}{Theorem}[section]
\newtheorem{Corollary}[Theorem]{Corollary}
\newtheorem{Lemma}[Theorem]{Lemma}
\newtheorem{Proposition}[Theorem]{Proposition}
 { \theoremstyle{definition}
\newtheorem{Definition}[Theorem]{Definition}

\newtheorem{Remark}[Theorem]{Remark} }

\newcommand{\cqfd}
{%
	\mbox{}%
	\nolinebreak%
	\hfill%
	\rule{2mm}{2mm}%
	\medbreak%
	\par%
}
\title{\bf\large Helix curves of the unit tangent bundle of a pseudo-Riemannian surface}
\author{\small Mohamed Tahar Kadaoui ABBASSI and Khadija BOULAGOUAZ}
\date{}
\begin{document}
	\maketitle
	

\textbf{Abstract}. In this paper, we classify helix (spacelike, timelike and null) curves, directed by the geodesic flow vector field, on the (3-dimensional) unit tangent bundle of a pseudo-Riemannian surface of constant Gaussian curvature endowed with a pseudo-Riemannian $g$-natural metric of Kaluza-Klein type. We find, in particular, that every such helix curve is, in fact, a circular helix in the senses that its curvature and torsion are constant.\\[2ex]


\vspace{0.6cm}

		
		\section*{Introduction and main results}
		

		In classical geometry, a curve $\lambda:I\rightarrow \mathbb{R}^3$, with unit speed, is a general helix if there exists a constant unit vector $V$ (the direction or axis of the helix) such that the scalar product $\left\langle \lambda'(t),V\right\rangle =\cos \theta$, for some constant angle $\theta$. Lancret's classical result states that if the torsion $\tau$ is non-zero, the ratio of curvature $\kappa$ to torsion remains constant along the curve. Since Lancret's time, helix curves have attracted significant attention due to their mathematical interest and applications in fields such as kinematics, DNA structure, carbon nanotubes, and liquid crystals. 
		
		The concept of helices has been extended to higher-dimensional and more general spaces by generalizing either the notions of angle and direction or Lancret's theorem. One of the earliest efforts was made by Hayden, who generalized helices in the context of $n$-dimensional Riemannian manifolds. He defined a helix curve as one parameterized by arc length with constant ratios between successive curvatures, characterizing these curves as making a constant angle with an arbitrary parallel vector field (the direction) along the curve (cf. \cite{hay1}). Barros, in \cite{Bar}, proposed a more general definition by requiring the helix direction to be a Killing vector field, rather than parallel. Barros also extended Lancret's theorem to Riemannian space forms, showing that in both spherical and Euclidean cases, general helices with non-zero torsion have curvature of the form $\kappa=a\tau\pm 1$ for some constant $a$,  while for the hyperbolic case, every general helix curve is a circular helix, that is both $\kappa$ and $\tau$ are constant.
		
		Later, the concept of a helix curve was extended to the realm of pseudo-Riemannian geometry. In such a manifold, there are three distinct types of curves: timelike, spacelike, and null curves. Timelike and spacelike curves behave similarly to their counterparts in Riemannian geometry. However, null curves require a different approach, as their arc length vanishes, making it impossible to normalize the tangent vector in the usual way. This absence of a standard parameter introduces challenges in constructing a Frenet frame for null curves, as it becomes dependent on the choice of both the screen vector bundle and the selected parameter. To address this issue, Ferr\'{a}ndez and colleagues introduced a new parameter, the pseudo-arc parameter, which normalizes the curve's acceleration rather than its speed. They developed an analogous Frenet frame, along with corresponding Frenet equations and curvatures, for null curves. These curvatures, referred to as the Cartan frame, Cartan equations, and Cartan curvatures of the null curve, remain invariant under Lorentzian transformations (cf. \cite{Fer_Gim}) Ferr\'{a}ndez also defined a Cartan curve as a null curve satisfying the Cartan equations, and a Cartan helix curve as a Cartan curve with constant curvatures.
		
		Additionally, since the concept of angle does not apply in a pseudo-Riemannian manifold, this challenge can be resolved by considering the scalar product $\left\langle \lambda',V\right\rangle$, where $V$ is a special vector field of constant length. This approach, adopted in \cite{Fer_Gim}, was used to define a generalized helix Cartan curve in Lorentz-Minkowski spaces and establish an analogue of Lancret's theorem for null curves. Specifically, a null Cartan curve is a generalized helix if, and only if, the first curvature function $\kappa_1$ and the ratios $\frac{\kappa_{2i+1}}{\kappa_{2i}}$ remain constant, where $\kappa_i$ represents the i-th curvatures of $\lambda$. In the 3-dimensional case, the definition and Lancret's theorem take on a simpler form: a null Cartan curve $\lambda:I\rightarrow \mathbb{R}^3_1$ is a generalized helix if there exists a non-zero constant vector $V$ such that the scalar product $\left\langle \lambda^\prime,V \right\rangle$ is constant along $\lambda$, which is characterized by the constancy of its lightlike curvature $\kappa_1$ \cite{Fer_Gim}.
		
		The study of helix curves becomes particularly significant when the ambient space contains a special vector field that can serve as the direction of the helix curves. This occurs, for instance, in contact or paracontact manifolds, where the Reeb vector field serves as the special vector field. In such cases, helices directed by the Reeb vector field correspond to the so-called slant curves with respect to the contact (or paracontact) structure. 

In this paper, the ambient space is the unit tangent bundle of a pseudo-Riemannian manifold endowed with an arbitrary pseudo-Riemannian metric induced by some pseudo Riemannian $g$-natural metric on the tangent bundle (cf. \cite{Abb-Sar}, \cite{Abb-Bou-Cal}). These metrics have a simpler form compared to those on general tangent bundles. Specifically, a metric $\tilde{G}$ is determined by four constants $a$, $b$, $c$ and $d$ in the following way
		\begin{equation} \label{g_natural_metric}
			\left\lbrace \begin{array}{l}
			\tilde{G}_{(x,u)}(X^h,Y^h) = (a + c)g_x(X, Y ) + dg_x(X, u)g_x(Y, u),\\
			\tilde{G}_{(x,u)}(X^h,Y^v) = bg_x(X, Y ),\\
			\tilde{G}_{(x,u)}(X^v, Y^v) = ag_x(X, Y ),
		\end{array}\right.
		\end{equation}
		for all $x \in M$ and $u, X, Y \in M_x$, where $X^h$ and $X^v$ stand for horizontal and vertical lifts of $X$, respectively.

If $b=0$ and $4a(a+c)=|a+c+d|$ (resp. $4a(a+c)=-|a+c+d|$) then we have a family of \textit{natural} contact pseudo-metric (resp. ($\varepsilon$)-paracontact metric) structures on the unit tangent bundle of a pseudo-Riemannian manifold (cf \cite{Abb_Cal2} (resp. \cite{Cal_15})). Such a structure is $K$-contact (resp. $K$-paracontact) if and only if the base manifold has constant positive (resp. negative) sectional curvature. Additionally, the Reeb vector field for these structures is collinear with the geodesic vector field on the unit tangent bundle, scaled by a constant factor. Consequently, helix curves on the unit tangent bundle directed along the geodesic vector field correspond to slant submanifolds with respect to one of these natural contact (resp. paracontact) metric structures.
		
		This paper focuses on helix curves in the unit tangent bundle directed by the geodesic vector field. This leads us to consider the natural contact pseudo-metric (resp. $(\varepsilon)$-paracontact metric) structures on the unit tangent bundle and to investigate the corresponding slant curves. For both geometric and technical reasons, we restrict our attention to a pseudo-Riemannian surface $M^2(\kappa)$ with non-zero constant Gauss curvature $\kappa$, and to $g$-natural metrics of the Kaluza-Klein type on its unit tangent bundle (i.e., metrics where the horizontal and vertical distributions are orthogonal, with $b=0$).

		To explore helix curves, we consider vertical, horizontal, and oblique curves. In the vertical case, we will prove that any smooth vertical curve in the unit tangent bundle of a pseudo-Riemannian surface, endowed with a pseudo-Riemannian Kaluza-Klein type metric, qualifies as a helix curve directed along the geodesic vector field.
		
		Next, we shift our focus to non-vertical smooth curves, beginning with the characterization of geodesic helices. Specifically, we show that geodesic helix curves correspond to parallel vector fields along geodesics, and that, in general, the metric on $T_1M^2(\kappa)$ is a Kaluza-Klein metric (i.e. $b=0$ and $d=0$). More precisely, we have the following results

		\begin{Theorem}\label{Geod-hel}
			Let $M^2(\kappa)$ be a pseudo-Riemannian surface with constant Gaussian curvature $\kappa \neq 0$, and let $T_1M^2(\kappa)$ denote its unit tangent bundle, equipped with a pseudo-Riemannian g-natural metric of Kaluza-Klein type (i.e. $b=0$). Consider a non-vertical smooth geodesic $\lambda=(x,V)$ on $T_1M^2(\kappa)$.
			\begin{enumerate}
				\item If $\lambda$ is timelike or spacelike and parameterized by its arc length, then $\lambda$ is a helix curve in $T_1M^2(\kappa)$ with the geodesic vector field as its direction if, and only if, $x$ is a non-constant timelike or spacelike geodesic, $V$ is a parallel vector field along $x$, and one of the following conditions holds:
				\begin{enumerate}
					\item $\lambda$ is a Legendre curve,
					\item $G$ is a Kaluza-Klein metric, i.e. $d=0$.
				\end{enumerate}
				\item If $\lambda$ is lightlike, then  $\lambda$ is a helix curve in $T_1M^2(\kappa)$ with the geodesic vector field as its direction if, and only if, $G$ is a Kaluza-Klein metric, $x$ is a non-constant lightlike geodesic, and $V$ is a parallel vector field along $x$.
			\end{enumerate}
		\end{Theorem}
		
		As an example, consider $M$ to be the Lorentzian pseudo-sphere $\mathbb{S}^2_1$, equipped with the restriction of the usual metric on $\mathbb{R}^3_1$, and $T_1M$ its unit tangent bundle, endowed with a pseudo-Riemannian $g$-natural metric of Kaluza-Klein type, with $c=0$. Let $f$ be the parameterization of $\mathbb{S}^2_1$ defined by 
$f(\theta,\phi)=(\sinh \theta,\cosh \theta\sin \phi,\cosh \theta \cos \phi)$ and $\delta\in \mathbb{R}$.
		Then the family of curve  $x(t)=(\delta t,1,\delta t),$
		 are lightlike geodesics of $M$. Therefore, if $d=0$ then the curves  $\lambda(t)=(x(t),V(t))$, where $V$ is the parallel transport of a unit vector $v\in T_{x(0)}M$ along $x$ are lightlike helix geodesics in $T_1M$ and if $d\neq 0$, $T_1M$ does not admit any lightlike helix curve.

		On the other hand, the family of curves  
			$x(t)=(\sqrt{1+\delta^2}\sinh t,\cosh t,\delta\sinh t),$
		where $\delta\in \mathbb{R}$, are timelike geodesics, and the curves	$\bar{x}(t)=(\sqrt{\varrho^2-1}\sin t,\cos t, \varrho \sin t)$,  
		 for $|\varrho|>1$, are spacelike geodesics. Therefore,
		\begin{itemize}
			\item if $d\neq 0$, let $v_0=(\delta,0,\sqrt{1+\delta^2})$ and $\bar{v}_0=(\varrho,0,\sqrt{\varrho^2-1}) $. Then the curves $$\lambda(t)=(x(t),V_0(t)) \quad \textup{and} \quad\bar{\lambda}(t)=(\bar{x}(t),\bar{V}_0(t)),$$ 
where $V_0$ (resp. $\bar{V}_1$) is the parallel transport along $x$ (resp. $\bar{x}$) of $v_0$ (resp. $\bar{v}_0$), are non-degenerate helix geodesics on $T_1M$,
			\item if $d=0$, then the curves 
$\lambda(t)=(x(t),V(t))$ and  $\bar{\lambda}(t)=(\bar{x}(t),\bar{V}(t)),$
where  $V$ (resp. $\bar{V}$) is the parallel transport of a unit vector $v\in T_{x(0)}M$ (resp. $\bar{v}\in T_{\bar{x}(0)}M$) along $x$ (resp. $\bar{x}$), are also non-degenerate helix geodesics on $T_1M$,.
		\end{itemize} 
	 In particular, if we take $x(t)=(\sinh(t),\cosh(t),0)$, then $x(t)$ and $-x(t)$ are timelike geodesics in $M$, consequently, for any $g$-natural Kaluza-Klein metric, the curve $\lambda(t)=(x(t),\dot{x}(t))$ and $-\lambda(t)$ as represented in Fig. \ref{fig:geod-hel-ex} below, are trivial helices in $T_1M$.
		\begin{figure}[h!]
			\begin{minipage}[c]{.46\linewidth}
				\centering
				\includegraphics[scale=0.6]{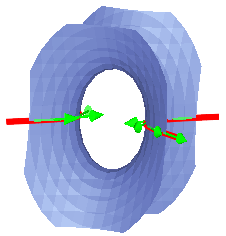}
			\end{minipage}
			\hfill
			\begin{minipage}[c]{.46\linewidth}
				\includegraphics[scale=0.6]{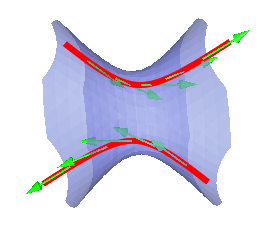}
			\end{minipage}
			\caption{Example of geodesic helices in $T_1\mathbb{S}^2_1$}
            \label{fig:geod-hel-ex}
		\end{figure}

		Unlike geodesic helices, we show that helix Frenet curves in the unit tangent bundle of a surface exist only when the metric on the unit tangent bundle is a non-Kaluza-Klein metric. In this case, every non-degenerate helix curve has constant curvature and torsion, while every lightlike helix curve has constant lightlike curvature. However, the precise classification of such helix curves depends on both the type of curve (whether it is horizontal or oblique) and its extrinsic geometric properties, such as torsion. 
		
		Let $(M^2(\kappa),g)$ be a pseudo-Riemannian surface of constant curvature $\kappa \neq 0$ and let $T_1M^2(\kappa)$ be its unit tangent bundle endowed with a pseudo-Riemannian $g$-natural metric $\tilde{G}$ of Kaluza-Klein type such that $\frac{a+c}{a}=\kappa$. We denote $\alpha=a(a+c)$ and $\varphi=a+c+d$.
		
		The helix horizontal Frenet curves with zero torsion are characterized as follows:
	
		\begin{Theorem}\label{th-hor-hel0}
			 Let $\lambda=(x,V)$ be a horizontal Frenet curve on $T_1M^2(\kappa)$ with zero torsion. Then $\lambda$ is a helix curve with direction the geodesic vector field on $T_1M^2(\kappa)$ if and only if the following assertions are satisfied
			\begin{enumerate}
				\item $d\neq 0$, $a+c-d \neq 0$;
				\item $x$ is a geodesic curve;  
				\item $g(\dot{x},V)^2=\frac{\varepsilon_\lambda}{2d}$; 
				\item the curvature $\kappa_\lambda$ of $\lambda$ is constant equals $\kappa_\lambda=\pm\sqrt{\frac{d-(a+c)}{4\varepsilon_1\alpha}}$;
			\end{enumerate}
			where  $\varepsilon_\lambda$ and $\varepsilon_1$ refer, respectively, to the first and second causal character of $\lambda$.
		\end{Theorem}
		
		Helix horizontal Frenet curves of non zero torsion are characterized as follows:
		
		\begin{Theorem}\label{th-hor-hel}
		 Let $\lambda=(x,V)$ be a horizontal Frenet curve of non zero torsion on $T_1M^2(\kappa)$. Then $\lambda$ is a helix curve with direction the geodesic vector field on $T_1M^2(\kappa)$ if and only if the following assertions are satisfied
			\begin{enumerate}
				\item $d\neq 0$, $4\alpha=\varphi$;
				\item $\varepsilon_2(\varepsilon -\varepsilon_\lambda\vartheta^2) >0$;  
				\item $\varepsilon_\lambda\varepsilon_1=\varepsilon_2$ and $x$ is a non-degenerate geodesic curve;
				\item both the curvature $\kappa_\lambda$ and the torsion $\tau$ of $\lambda$ are constants satisfying
				\begin{equation*}\label{Curv-tors-hor-frenet}
					\kappa_\lambda=\pm\frac{2\vartheta d}{\varphi}\sqrt{\varepsilon_2 (\varepsilon-\varepsilon_\lambda\vartheta^2)},
					\quad \tau=\pm\left(\varepsilon\varepsilon_\lambda-\frac{2d\vartheta^2}{\varphi}\right);
				\end{equation*}
			\end{enumerate}
			where $\varepsilon=\frac{|\varphi|}{\varphi}$, $\varepsilon_\lambda, \varepsilon_1, \varepsilon_2$ refer, respectively, to the first, second and third causal characters of $\lambda$ and $\vartheta=\varepsilon\sqrt{|\varphi|}g(\dot{x},V)$.
		\end{Theorem}

On the other hand, helix oblique Frenet curves with zero torsion, along constant speed curves, are characterized as follows:		
		\begin{Theorem}\label{hel-obl-fren0}
			Let $\lambda=(x,V)$ be a timelike or a spacelike oblique Frenet curve on $T_1M^2(\kappa)$ of zero torsion, such that the curve $x$ has a constant speed, and let $\sigma:=g(\dot{x},\dot{x})$. Then $\lambda$ is a helix with direction the geodesic vector field if, and only if, the following assertions are satisfied
			\begin{enumerate}
				\item $d\neq 0$;
				\item $x$ is a space-like (pseudo-)Riemannian circle and $\sigma=\frac{\varepsilon\vartheta^2}{\varphi}$;
				\item $V=\pm\frac{1}{\sqrt{\sigma}}\dot{x}$;
				\item the curvature $\kappa_\lambda$ of $\lambda$ is constant equals $\kappa_\lambda=\pm\left(-\frac{ d(a+c-d)}{2\varepsilon_1\varepsilon_\lambda\alpha|\varphi|}\right)^{1/2};$
			\end{enumerate}	
				where $\varepsilon=\frac{|\varphi|}{\varphi}$, $\varepsilon_\lambda$ and $\varepsilon_1$ refer, respectively, to the first and second causal characters of $\lambda$ and $\vartheta=\varepsilon\sqrt{|\varphi|}g(\dot{x},V)$.
		\end{Theorem}
		
		As example, let $M=\mathbb{S}^2$ represent the Euclidean sphere in $\mathbb{R}^3$. The curve $x(t)=\frac{1}{\sqrt{2}}(\cos(t),1,\sin(t))$ 
is a Riemannian circle with $\sigma=\frac{1}{2}$. We equip $T_1M$  with a metric \eqref{g_natural_metric}  satisfying $a=a+c=1$, $d=3$. In this setting, the curve $\lambda(t)=(x(t),\sqrt{2}\dot{x}(t))$ is an oblique helix in $(T_1M,\tilde{G})$ with an angle $\vartheta=\sqrt{2}$ (cf. Fig. \ref{fig:obl-hel-ex} below). Additionally, its curvature is $\kappa_\lambda=\frac{\sqrt{3}}{2}$.
		
		\begin{figure}[h!]
			\begin{minipage}[c]{.46\linewidth}
				\centering
				\includegraphics[scale=0.4]{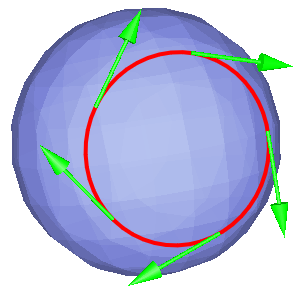}
			\end{minipage}
			\hfill
			\begin{minipage}[c]{.46\linewidth}
				\includegraphics[scale=0.4]{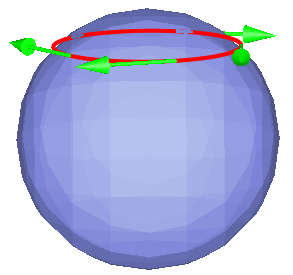}
			\end{minipage} 
			\caption{Example of an oblique helix curve in $T_1\mathbb{S}^2$.}
            \label{fig:obl-hel-ex}
		\end{figure}

		Finally, helix oblique Frenet curves with non vanishing torsion, along constant speed curves, are characterized as follows:
		
		\begin{Theorem}\label{hel-obl-fren}
			Let $\lambda=(x,V)$ be a timelike or a spacelike oblique Frenet curve on $T_1M$ of non zero torsion, such that the curve $x$ has a constant speed, and let $\sigma:=g(\dot{x},\dot{x})$. Then $\lambda$ is a helix with direction the geodesic vector field if, and only if, the following assertions are satisfied
			\begin{enumerate}
				\item $d\neq 0$, $\varepsilon_1\varepsilon_2\alpha >0$;
				\item $x$ is a space-like (pseudo-)Riemannian circle and $\sigma=\frac{\varepsilon\vartheta^2}{\varphi}$;
				\item $\varepsilon\varepsilon_\lambda=1$, $\varepsilon_2(\varepsilon-\varepsilon_\lambda\vartheta^2)>0$;
				\item $V=\pm\frac{1}{\sqrt{\sigma}}\dot{x}$;
				\item both the curvature $\kappa_\lambda$ and the torsion $\tau$ of $\lambda$ are constants satisfying 
				\begin{equation*}
					\kappa_\lambda=\pm d\vartheta\sqrt{\frac{\varepsilon_1(1-\vartheta^2)}{\alpha\varphi}},\qquad 				
					\tau=\pm\frac{\varepsilon\varphi-2\varepsilon_\lambda d\vartheta^2}{2\sqrt{|\alpha\varphi|}};
				\end{equation*}
				where  $\varepsilon=\frac{|\varphi|}{\varphi}$, $\varepsilon_\lambda, \varepsilon_1, \varepsilon_2$ refer, respectively, to the first, second and third causal characters of $\lambda$ and $\vartheta=\varepsilon\sqrt{|\varphi|}g(\dot{x},V)$.
			\end{enumerate}	
		\end{Theorem}
		
In the lightlike case, non-geodesic helix horizontal and oblique curves are described by the following theorems:
		\begin{Theorem}\label{null-hor-hel}
			Let $\lambda=(x,V)$ be a Cartan Frenet lightlike horizontal curve on $T_1M^2(\kappa)$, parameterized by the pseudo-arc parameter. Then $\lambda$ is helix curve with direction the geodesic vector field on $T_1M^2(\kappa)$ if and only if 
			\begin{enumerate}
				\item $d\neq 0$ and $\varphi=-4\alpha$; 
				\item $g(\dot{x},V)=\frac{1}{2|d|}$;
				\item either $x$ is a non-degenerate geodesic or $a+c=3d$ and $\ddot{x}$ is a lightlike vector field along $x$;
				\item the lightlike curvature $\kappa_\lambda$ of $\lambda$ is equal to $\kappa_\lambda = -\frac{d(a+c)}{|d\varphi|}$.
			\end{enumerate}
		\end{Theorem}
		
		\begin{Theorem}\label{null-obl-hel}
			Let $\lambda=(x,V)$ be an oblique Cartan Frenet lightlike curve parameterized by the pseudo-arc parameter on $T_1M^2(\kappa)$, such that the curve $x$ has a constant speed, and let $\sigma:=g(\dot{x},\dot{x})$. $\lambda$ is a helix curve with direction the geodesic vector field if and only if the following assertions hold
			\begin{enumerate}
				\item $d=a+c\neq 0$;  
				\item $\sigma=\frac{\varepsilon}{2d}$,  and $V=\pm\sqrt{2\varepsilon d}\dot{x}$; 
				\item $x$ is a pseudo-Riemannian circle; 
				\item the lightlike curvature $\kappa_\lambda$ of $\lambda$ is equal to $\kappa_\lambda=\frac{\varepsilon }{2}$.
			\end{enumerate} 
		\end{Theorem}
		

		\section{Preliminaries}
		
		
		\subsection{contact and paracontact structures on a manifold}
		
		
		Let us now review some basic information about contact and paracontact structures on a manifold. 
		
		An \textit{almost contact structure} (resp. \textit{almost paracontact structure}) on an odd-dimensional smooth manifold $\overline{M}$ consists of a triplet $(\phi,\xi,\eta)$, where $\eta$ is 1-form, $\xi$  is a global vector field and $\phi$ is a $(1,1)$-tensor field that satisfies the following conditions:
		\begin{enumerate}
			\item $\eta(\xi)=1$,
			\item $\phi^2=-Id+\eta\otimes\xi$ (resp. $\phi^2=Id-\eta\otimes\xi$).
		\end{enumerate}
		It follows from these conditions that $\phi\xi=0$ and $\eta\circ \phi=0$.
		
		A pseudo-Riemannian metric $\bar{g}$ on $\overline{M}$ is said to be \textit{compatible with the almost contact structure} (resp. \textit{paracontact structure}) $(\phi,\xi,\eta)$ if 
		\begin{equation}\label{cont-para}
			\bar{g}(\phi X,\phi Y)=\bar{g}(X,Y)-\varepsilon\eta(X)\eta(Y)\quad \textup{(resp.} \quad \bar{g}(\phi X,\phi Y)=-\bar{g}(X,Y)+\varepsilon\eta(X)\eta(Y)\mathrm{)},
		\end{equation}
		where $\varepsilon=\pm 1=\bar{g}(\xi,\xi)$. 

A smooth manifold equipped with an almost contact (or paracontact) structure and a compatible pseudo-Riemannian metric is called an \textit{almost contact pseudo-metric manifold} (resp. \textit{($\varepsilon$)-almost paracontact metric manifold}). If  the compatible pseudo-Riemannian metric further satisfies 
		\begin{equation}\label{contact-eta}
			\bar{g}(X,\phi Y)=(d\eta)(X,Y),
		\end{equation}
		then $\eta$ is a \textit{contact form} (resp. \textit{paracontact form}) on $\overline{M}$, $\xi$ is \textit{the associated Reeb vector field}, $\bar{g}$ is an \textit{associated metric} and $(\overline{M},\eta,\bar{g})$ is called  a \textit{contact pseudo-metric} manifold (resp. \textit{($\varepsilon$)-paracontact metric}). The structure $(\overline{M},\eta,\bar{g})$ is referred to as a \textit{$K$-contact} (resp. \textit{$K$-paracontact}) manifold if $\xi$ is a Killing vector field. 

In \cite{Cal_Per} the authors showed that $(\overline{M},\eta,\bar{g})$ being $K$-contact  is equivalent to the vanishing of the tensor field $h:=\frac{1}{2}\mathcal{L}_\xi \varphi$. Furthermore, similar to the Riemannian case, we have
		\begin{equation}\label{contact-sasaki}
			\nabla_\xi\varphi=0 \quad \textup{ and} \quad \nabla_X\xi=-\varepsilon\phi X-\phi hX \quad \textup{(resp.} \quad \nabla_X\xi=-\varepsilon\phi X+\phi hX\mathrm{)}.
		\end{equation}
		Additionally, $h$ is self-adjoint, meaning $h\phi=-\phi h$ and $h\xi=0$. 

		A manifold $(\overline{M},\eta,\bar{g})$ is called \textit{Sasakian} (resp. \textit{parasasakian}) if $[\phi,\phi]+2\eta\otimes\xi=0$, where $[\phi,\phi]$ is defined by 
		$$[\phi,\phi](X,Y)=\phi^2[X,Y]+[\phi X,\phi Y]-\phi[\phi X,Y]-\phi[X,\phi Y],$$
		for any vector fields $X,Y\in \mathfrak{X}(M)$.  
		\begin{Proposition}(\cite{Cal_Per})
			Let $(\overline{M},\eta,\bar{g})$ be an almost contact pseudo-metric manifold. The manifold $(\overline{M},\eta,\bar{g})$ is Sasakian (resp. parasasakian) if, and only if, 
			$$(\nabla_X\varphi)Y=\bar{g}(X,Y)\xi-\varepsilon\eta(Y)Z \quad \mathrm{(resp.}\quad (\nabla_X\varphi)Y=-\bar{g}(X,Y)\xi+\varepsilon\eta(Y)Z\:\mathrm{)}, $$
			for all $X, Y \in \mathfrak{X}(M)$
		\end{Proposition}
		
		
		\subsection{Riemannian geometry of tangent bundles}
		
		
		Let $\nabla$ be the Levi-Civita connection of $g$. Then, at any point $(x,u)\in TM$, the tangent space of $TM$ splits into horizontal and vertical subspaces with respect to $\nabla$:
		$$(TM)_{(x,u)}=H_{(x,u)}\oplus V_{(x,u)}.$$    \par
		
		For any vector $X\in M_x$, there exists a unique vector $X^h  \in H_{(x,u)}$ such that $p_* X^h =X$. We call $X^h$ the \emph{horizontal lift} of $X$ to the point $(x,u)\in TM$. The \emph{vertical lift} of a vector $X\in M_x$ to $(x,u)\in TM$ is a vector $X^v  \in V_{(x,u)}$ such that $X^v (df) =Xf$, for all functions $f$ on $M$. Here, we consider the $1$-forms $df$ on $M$ as functions on $TM$, $(df)(x,u)=uf$. 

Note that the map $X \to X^h$ is an isomorphism between $M_x$ and $H_{(x,u)}$. Similarly, the map $X \to X^v$ is an isomorphism between $M_x$ and $V_{(x,u)}$. Therefore, each tangent vector $\tilde Z \in (TM)_{(x,u)}$ can be uniquely expressed as $\tilde Z =X^h + Y^v$, where $X,Y \in M_x$.
		
		On the tangent bundle of a pseudo-Riemannian manifold, we can construct a large class of metrics called $g$-natural metrics, which are determined as follows:
		
		\begin{Proposition}\label{$g$-nat} \cite{Abb-Bou-Cal}
			Let $(M,g)$ be a pseudo-Riemannian manifold, and let $G$ be a $g$-natural metric on $TM$. Then there exist functions $\alpha_i$, $\beta_i:\mathbb{R} \rightarrow \mathbb{R}$, $i=1,2,3$, such that for every $u$, $X$, $Y\in M_x$, the metric $G$ is given by:
			\arraycolsep1.5pt
			\begin{equation}\label{g_exp}
				\left\lbrace
				\begin{array}{rcl}
					G_{(x,u)}(X^h,Y^h)& = & (\alpha_1+ \alpha_3)(\rho) g_x(X,Y)
					+ (\beta_1+ \beta_3)(\rho)g_x(X,u)g_x(Y,u),\\
					G_{(x,u)}(X^h,Y^v)& = & \alpha_2 (\rho) g_x(X,Y)
					+  \beta_2 (\rho) g_x(X,u)g_x(Y,u), \\
					G_{(x,u)}(X^v,Y^v)& = & \alpha_1 (\rho) g_x(X,Y)
					+ \beta_1 (\rho) g_x(X,u)g_x(Y,u),
				\end{array}
				\arraycolsep5pt \right.
			\end{equation}
			where $\rho =g_x(u,u)$.
		\end{Proposition}
		
		Pseudo-Riemannian $g$-natural metrics are characterized by the following proposition:
		
		\begin{Proposition}\label{riem-nat} \cite{Abb-Bou-Cal}
			A $g$-natural metric $G$ on the tangent bundle of a pseudo-Riemannian manifold $(M,g)$, defined by the functions in Proposition \ref{$g$-nat} is non-degenerate if and only if
		$\alpha(t) \neq 0$, $\phi(t) \neq 0$ for all $ t \in \mathbb{R}, $
			where the functions $\phi_i$, $\alpha$, $\phi$ are defined by
\begin{itemize}
  \item $\phi_i(t)=\alpha_i(t) +t\beta_i(t)$, for all $t \in \mathbb{R}$,
  \item $\alpha=\alpha_1(\alpha_1+\alpha_3)-\alpha_2^2$;
  \item $\phi=\phi_1(\phi_1+\phi_3)-\phi_2^2$.
\end{itemize}
		\end{Proposition}
		
		The wide class of $g$-natural metrics includes several well-known metrics (both Riemannian and non-Riemannian) on $TM$. Specifically:
		\begin{itemize}
			\item the {\em Sasaki metric} $g_S$ is obtained for $\alpha _1 =1$ and $\alpha _2 = \alpha _3 = \beta _1 =\beta _2 = \beta _3 =0$.
			\item {\em Kaluza--Klein metrics}, commonly defined on principal bundles \cite{Wood1} {(see also \cite{Loubeau})}, are obtained for
			$\alpha _2 = \beta _2 = \beta _1 +\beta _3 = 0$.
			\item {\em Metrics of Kaluza--Klein type} are defined by the geometric condition of orthogonality between horizontal and vertical distributions. Thus, a $g$-natural metric $G$ is of Kaluza-Klein type if $\alpha _2=\beta _2 =0$.
		\end{itemize}
		
		
		\subsection{$g$-natural metrics and associated contact and paracontact structures on unit tangent bundles of pseudo-Riemannian manifolds} 
		
		
		\label{sect-2}

		Let $(M,g)$ be a pseudo-Riemannian manifold with signature $(n-k,k)$. The unit tangent sphere bundle over $M$, denoted by $T_1M$, is defined as the hypersurface of the tangent bundle $TM$ by the condition 
$$T_1M=\left\lbrace (x,u)\in TM:g_x(u,u)=1\right\rbrace.$$ 
The tangent space of $T_1M$ at a point $(x,u)\in T_1M$ is given by 
$$T_{(x,u)}T_1M =\left\lbrace X^h + Y^v: X, Y \in T_xM,\; g(Y, u) = 0	\right\rbrace.$$
		
Let $G$ be a $g$-natural metric on $TM$, as defined by the functions in Proposition \ref{$g$-nat}, and let $\tilde{G}$ denote its restriction to $T_1M$. Following the same approach as in \cite{Abb_Kow0} for Riemannian base manifolds, we can characterize $g$-natural metrics on the unit tangent bundle of pseudo-Riemannian manifolds. It is straightforward to see that $\tilde{G}$ is fully determined by equation \eqref{g_natural_metric}: 
		$$\left\lbrace \begin{array}{l}
			\tilde{G}_{(x,u)}(X^h,Y^h) = (a + c)g_x(X, Y ) + dg_x(X, u)g_x(Y, u),\\
			\tilde{G}_{(x,u)}(X^h,Y^v) = bg_x(X, Y ),\\
			\tilde{G}_{(x,u)}(X^v, Y^v) = ag_x(X, Y ),
		\end{array}\right.$$
		where $a = \alpha_1(1)$, $c = \alpha_3(1)$, $d = \beta_3(1)$ and $b = \alpha_2(1)$.\\
		
		In particular, 
		\begin{itemize}
			\item {\em Kaluza-Klein metrics} are obtained when $b=d=0$;
			\item {\em metrics of Kaluza-Klein type} are defined when $b=0$.
		\end{itemize}
		
		We denote $\varphi=a+c+d$ and $\alpha=a(a+c)-b^2$. The vector field, defined for all $(x,u) \in T_1M$, by 
		$$N^G_{(x,u)}=\frac{1}{\sqrt{|\varphi\phi|}}\left[-bu^h+\varphi u^v\right],$$ 
		is unitary and normal to $T_1M$ at any point $(x,u)\in T_1M$. The \textit{tangential lift} of a vector field $X\in\mathfrak{X}(M)$ at $(x,u)$ is then defined by    
		\begin{equation}\label{tang-lift}
			X^{t}=X^v-\frac{1}{\varphi}g(X,u)(-bu^h+\varphi u^v).
		\end{equation}
		
		As a result, the tangent space to $T_1M$ at $(x,u)$ can be expressed as (\cite{Cal_15})
		$$T_{(x,u)}T_1M =\left\lbrace X^h + Y^{t}: X, Y \in T_xM, \; g(Y, u) = 0	\right\rbrace.$$
		\begin{Proposition}
			The metric $\tilde{G}$ is non degenerate on $T_1M$ if and only if $\alpha\neq 0$ and $\varphi\neq 0$.
		\end{Proposition}
		
		\begin{Proposition}\label{signature}
			The signature of a $g$-natural metric on $T_1M$ is determined as follows
			\begin{enumerate}
				\item 	$(2n-2k-1, 2k)$ if $\alpha>  0$,  $a+c>0$ and $a>0$,
				\item $(2k+1,2n-2k-2)$ if $\alpha>  0$, $\varphi >0$, $a+c<0$ and $a<0$,
				\item $(2k,2n-2k-1)$ if $\alpha>0$, $\varphi <0$, $a+c<0$ and $a<0$,
				\item $(n,n-1)$ if $\alpha< 0$ and $\varphi>0$,
				\item $(n-1,n)$ if $\alpha< 0$ and $\varphi<0$.
			\end{enumerate}	
		\end{Proposition}

		
		In this paper, we explore Sasakian (resp. parasasakian) contact pseudo-metric (resp. ($\varepsilon$)-paracontact metric) structures $(\tilde{G}, \tilde{\eta},\tilde{\phi},\tilde{\xi})$ on $T_1M$, where the compatible metrics $\tilde{G}$ are pseudo-Riemannian $g$-natural metrics.  To this end, we assume that the characteristic vector field $\tilde{\xi}$ is a unit tangent vector field on $T_1M$, collinear with the geodesic vector field at each point $(x,u)\in T_1M$. Depending on the causal character we wish to assign to $\tilde{\xi}$, we have
		\begin{equation}\label{xi}
			\tilde{\xi}_{(x,u)}=\frac{1}{\sqrt{\varepsilon\varphi}}u^h,
		\end{equation}
		for each $(x,u)\in T_1M$, ensuring that $\pm 1= \varepsilon=\tilde{G}(\tilde{\xi},\tilde{\xi})$. It follows then that $\varepsilon\varphi=|\varphi|$. 
		
		For the contact (resp. paracontact) metric case, the first author and G. Calvaruso \cite{Abb_Cal2} (resp. G. Calvaruso and V. Martin-Molina \cite{Cal_15}) proved that $(\tilde{G}, \tilde{\eta},\tilde{\phi},\tilde{\xi})$ is a contact metric structure if and only if $b=0$ and $\varphi=4\alpha$ (resp. $\varphi=-4\alpha$). Extending these arguments to our general case, we get the following proposition:
		\begin{Proposition}\label{contact}
			$(\tilde{G}, \tilde{\eta},\tilde{\phi},\tilde{\xi})$ is a contact pseudo-metric (resp. ($\varepsilon$)-paracontact metric) structure if and only if $b=0$ and $|\varphi|=4\alpha$ (resp. $|\varphi|=-4\alpha$).
		\end{Proposition} 
		
		Consequently, we will assume that $\tilde{G}$ is a Kaluza-Klein type metric. Following similar reasoning as in the contact (resp. paracontact) metric cases (cf. \cite{Abb_Cal2} (resp. \cite{Cal_15})), straightforward to show that the function $\tilde{\eta}$ and the tensor field $\tilde{\phi}$ on $T_1M$ are given by 
		\begin{equation}\label{eta}
			\tilde{\eta}(X^h)=\sqrt{|\varphi|}g(X,u),\quad \tilde{\eta}(X^t)=0,
		\end{equation}
		and \begin{equation}\label{phi}
			\tilde{\phi}(X^h)=\frac{\sqrt{|\varphi|}}{2a}X^t,\quad 
			\tilde{\phi}(X^t)=-\frac{\sqrt{|\varphi|}}{2(a+c)}\left[X^h-g(X,u)u^h\right],
		\end{equation}
		for any $X\in T_xM $. 
		
		\begin{Remark}
			As a consequence of Proposition \ref{contact}, if $(\tilde{G}, \tilde{\eta},\tilde{\phi},\tilde{\xi})$ is a contact pseudo-metric (resp. ($\varepsilon$)-paracontact metric) structure, then $\alpha>0$ (resp. $\alpha<0$). Thus, by Proposition \ref{signature}, we deduce: 
			\begin{itemize}
				\item [$\bullet$] If $(\tilde{G}, \tilde{\eta},\tilde{\phi},\tilde{\xi})$ is a contact pseudo-metric structure, then the signature of $\tilde{G}$ is one of the following $(2n-2k-1,2k)$, $(2k+1,2n-2k-2)$, $(2k,2n-2k-1)$, or $(2n-2k-2,2k+1)$. In the first two cases, $\tilde{\xi}$ is spacelike and, in the latter cases, it is timelike. 
				\item [$\bullet$] if $(\tilde{G}, \tilde{\eta},\tilde{\phi},\tilde{\xi})$ is a ($\varepsilon$)-paracontact metric structure, then the signature of $\tilde{G}$ is either $(n,n-1)$ or $(n-1,n)$. In the first case, $\tilde{\xi}$ is spacelike and in the second it is timelike.  
			\end{itemize}
		\end{Remark}
		
		In the same manner as in \cite{Abb_Cal2} (resp. \cite{Cal_15}) for the contact (resp. paracontact) metric cases, we can show that Sasakian (resp. parasasakian) contact pseudo-metric (resp. ($\varepsilon$)-paracontact metric) structures on $T_1M$ are characterized by the following
		\begin{Proposition}\label{K-contact}
			The contact pseudo-metric (resp. ($\varepsilon$)-paracontact metric) structure $(\tilde{\eta},\tilde{G},\tilde{\phi},\tilde{\xi})$ is $K$-contact (resp. $K$-paracontact) if and only if $(M,g)$ has a constant sectional curvature $\kappa=\frac{a+c}{a}>0$ (resp. $\kappa=\frac{a+c}{a}<0$). Moreover, in this case, $(\tilde{\eta},\tilde{G},\tilde{\phi},\tilde{\xi})$ is Sasakian (resp. parasasakian).
		\end{Proposition}
		
		Henceforth, we consider a pseudo-Riemannian manifold $(M(\kappa),g)$ of constant sectional curvature $\kappa$ and a pseudo-Riemannian $g$-natural metric $\tilde{G}$ of Kaluza-Klein type ($b=0$) on $T_1M(\kappa)$, satisfying $a+c=a\kappa$ and $|\varphi|=4|\alpha|$. Given that $b=0$, then by \eqref{tang-lift} the tangential lift of any vector field $X\in \mathfrak{X}(M)$ is  
		$$X^{t}=X^v-g(X,u) u^v.$$ 
		As a result, the metric $\tilde{G}$ is explicitly described  by 
		\begin{equation}
			\left\lbrace\begin{array}{l}
				\tilde{G}_{(x,u)}(X^h,Y^h) = (a + c)g_x(X, Y ) + dg_x(X, u)g_x(Y, u),\\
				\tilde{G}_{(x,u)}(X^h,Y^t) = 0,\\
				\tilde{G}_{(x,u)}(X^t, Y^t) = ag_x(X, Y )-ag_x(X,u)g_x(Y,u),
			\end{array}\right.
			\label{nat_met}
		\end{equation}
		
		Using similar methods as in \cite{Abb_Cal2} for a Riemannian base manifold, we prove that the Levi-Civita connection of the unit tangent bundle of a pseudo-Riemannian manifold, equipped with an arbitrary pseudo-Riemannian $g$-natural metric, follows the same characterization as in Proposition 5 from \cite{Abb_Cal2}. In our specific case, we have
		\begin{Proposition} 
			The Levi-Civita connection on $(T_1M(\kappa),\tilde{G})$ is totally described by: 
			\begin{align*}
					\tilde{\nabla}_{X^h}Y^h&=\left\lbrace\nabla_XY\right\rbrace^h+\left\lbrace -\frac{\varphi}{2a}g(Y,u)X+\frac{a+c-d}{2a}g(X,u)Y\right\rbrace^t,\\
					\tilde{\nabla}_{X^h}Y^t&=\left\lbrace \frac{d-(a+c)}{2(a+c)}g(X,u)Y+\frac{1}{2}g(X,Y)u-\frac{d}{2(a+c)}g(X,u)g(Y,u)u\right\rbrace^h+\left\lbrace \nabla_XY\right\rbrace^t,\\
					\tilde{\nabla}_{X^t}Y^h&=\left\lbrace \frac{d-(a+c)}{2(a+c)}g(Y,u)X+\frac{1}{2}g(X,Y)u-\frac{d}{2(a+c)}g(X,u)g(Y,u)u\right\rbrace^h,\\
					\tilde{\nabla}_{X^t}Y^t&=-\left\lbrace g(Y,u)X\right\rbrace^t.			
			\end{align*}	
			\label{cov-der}
		\end{Proposition}

		
		\section{Discussion and proofs of the main results}
		

		
		\subsection{Helices on $T_1M^2(\kappa)$}
		
		\label{hel-curves}
		\begin{Definition}\label{Helix-curve}
			Let $(M,g)$ be an $n$-dimensional pseudo-Riemannian manifold, $V$ a constant length vector field on $M$, and $\lambda:I \rightarrow M$ a smooth curve parametrized by its (pseudo-)arc length. We say that $\lambda$ is \emph{helix curve with direction $V$} if $g(\lambda^\prime,V)$ remains constant along $\lambda$.
		\end{Definition}

		In this section, our objective is to characterize helix curves on $(T_1M^2(\kappa),\tilde{G})$ directed by the geodesic vector field. In other words, by Definition \ref{Helix-curve}, these are curves $\lambda:I \rightarrow M^2(\kappa)$, parameterized by their (pseudo-)arc length, such that $G(\lambda^\prime,u^h)$ remains constant along $\lambda$.
		
		We begin by considering vertical curves. Let $\lambda=(x,V)$ be a vertical smooth curve in $T_1M^2(\kappa)$ parameterized by its arc-length. In this case, $x$ is a constant curve, equal to $x_0$, and $\lambda$ lies entirely within the fiber of $T_1M^2(\kappa)$ over $x_0$, which forms the (1-dimensional) unit circle of $T_{x_0}M^2(\kappa)$. Therfore $\lambda$ is locally the arc-length parametrization of this circle. Consequently $\lambda^\prime(t)$ is a vertical vector, for all $t$, and is then orthogonal to the geodesic vector field. Thus $\lambda$ is a Legendre curve, as a result, any vertical smooth curve in $T_1M^2(\kappa)$ is a helix curve with direction along the geodesic vector field. We now turn to the investigation of non-vertical smooth curves.
		
		Let $\lambda=(x,V)$  be a non-vertical smooth curve of $T_1M^2(\kappa)$, either null or parameterized by its (pseudo-)arc length. The velocity of this curve is given by:
		\begin{equation}\label{vel-0}
			\lambda'=\dot{x}^h+(\nabla_{\dot{x}}V)^t,
		\end{equation}
		such that
		\begin{equation}\label{vel-1}
			\varepsilon_\lambda= \tilde{G}(\lambda',\lambda')=(a+c)g(\dot{x},\dot{x})+dg(\dot{x},V)^2+ag(\dot{V},\dot{V}),
		\end{equation}
		where $\varepsilon_\lambda =\pm1$ or $0$, depending on the causal character of $\lambda$. 
		
		\begin{Proposition}\label{prop1}
			$\lambda$ is a helix curve in $T_1M^2(\kappa)$ with direction the geodesic vector field (or equivalently, a slant curve with respect to the Sasakian contact pseudo-metric (resp. parasasakian ($\varepsilon$)-paracontact metric) structure $(\tilde{\eta},\tilde{G},\tilde{\phi},\tilde{\xi})$, as given by \eqref{xi}, \eqref{eta} and \eqref{phi}, i.e. $\tilde{G}(\lambda',\tilde{\xi})$ is a constant $\vartheta$) if and only if $g(\dot{x},V)$ is constant.
		\end{Proposition}
		\emph{Proof of Proposition \ref{prop1}:}
			Since $
				\tilde{G}(\lambda',\tilde{\xi})=\frac{\varphi}{\sqrt{|\varphi|}}g(\dot{x},V)=\varepsilon\sqrt{|\varphi|}g(\dot{x},V),$
			then
			\begin{equation}
				g(\dot{x},V)=\frac{\varepsilon}{\sqrt{|\varphi|}}\vartheta,
				\label{eq4}
			\end{equation}
			which is constant if and only if $\vartheta$ is constant.
		\cqfd

		
		\subsection*{Geodesic helices on $T_1M^2(\kappa)$}
		

		\emph{Proof of Theorem \ref{Geod-hel}:}
		Suppose $\lambda$ is a geodesic. Using Proposition \ref{cov-der}, formula \eqref{eq4}, and $g(V,\dot{V})=0$ ($V$ being unitary), we derive the following system of equations:
		\begin{equation}
			\left\lbrace
			\begin{split}
				&	\ddot{x}-\frac{\varepsilon(a+c-d)}{(a+c)\sqrt{|\varphi|}}\vartheta\dot{V}+g(\dot{V},\dot{x})V=0,\\
				&	\ddot{V}-\frac{\varepsilon d}{a\sqrt{|\varphi|}}\vartheta\dot{x}=0.
			\end{split}
			\right.
			\label{eq_geo1}
		\end{equation} 
		Taking the scalar product of the second equation with $V$, we find $g(\ddot{V},V)=\frac {d\vartheta^2}{a|\varphi|}.$
		Since $g(\dot{V},V)=0$, it follows that  
		\begin{equation*}
			g(\dot{V},\dot{V})=-g(\ddot{V},V)=-\frac {d\vartheta^2}{a|\varphi|}.
		\end{equation*}
		Substituting this result, along with \eqref{eq4}, into \eqref{vel-1}, we find 
		$g(\dot{x},\dot{x})=\frac{\varepsilon_\lambda}{a+c},$
		which shows that the curve $x$ has a constant length. If $\lambda$ is lightlike, then $\varepsilon_\lambda=0$, and consequently $x$ is lightlike.
		Now, taking the scalar product of the first equation in \eqref{eq_geo1} by $\dot{x}$, and using \eqref{eq4}, we find
		\begin{equation*}
			0=g(\ddot{x},\dot{x})=-\frac{\varepsilon\vartheta d}{(a+c)\sqrt{|\varphi|}}g(\dot{V},\dot{x}).
		\end{equation*}
		This implies that either $\vartheta=0$, $d=0$ or $g(\dot{V},\dot{x})=0$.
		\begin{enumerate}
			\item If $\vartheta=0$, then equations (\ref{eq_geo1}) simplify to
			\begin{equation}
					\ddot{x}+g(\dot{V},\dot{x})V=0,\quad \textup{and} \quad
						\ddot{V}=0.
				\label{eq_geo2}
			\end{equation}
			In this case, we also have $g(\dot{V},\dot{V})=0$, which implies that either $\dot{V}=0$, or $\dot{V}$ is a lightlike vector field. 
			In the latter case, $\{V,\dot{V}\}$ forms a local orthogonal frame along $x$, and since $g(\dot{V},\dot{V})=g(\dot{V},V)=0$, 
			$g$ becomes degenerate, which contradicts our assumption. Thus, we conclude that $\dot{V}=0$.
			\begin{itemize}
				\item If $\lambda$ is spacelike or timelike, then equation (\ref{eq_geo2}) implies that $x$ is a geodesic.
				\item If $\lambda$ is lightlike, then $\lbrace \dot{x},V\rbrace$ forms a local frame along $x$ with $g(\dot{x},\dot{x})=g(\dot{x},V)=0$, 
				contradicting the non-degeneracy of $g$. Hence, $\vartheta\neq 0$.
			\end{itemize}
			\item If $d=0$, then $g(\dot{V},\dot{V})=0$, and following the same reasoning as in the previous case, we find $\dot{V}=0$. Therefore, equations (\ref{eq_geo1}) reduce to $\ddot{x}=0$, that is $x$ is a geodesic.
			\item Now, suppose $d\neq 0$ and $\vartheta\neq 0$. In this case, $g(\dot{V},\dot{x})=0$, and equations (\ref{eq_geo1}) become
			\begin{equation}
				\left\lbrace
				\begin{split}
					&	\ddot{x}-\frac{\varepsilon(a+c-d)}{(a+c)\sqrt{|\varphi|}}\vartheta\dot{V}=0,\\
					&	\ddot{V}-\frac{\varepsilon d}{a\sqrt{|\varphi|}}\vartheta\dot{x}=0.
				\end{split}
				\right.
				\label{eq_geo22}
			\end{equation} 
			Taking the scalar product of the first equation by $\dot{V}$ and the second by  we obtain $\dot{x}$
			\begin{equation}
				g(\ddot{x},\dot{V})=-\frac{\varepsilon d(a+c-d)}{\alpha|\varphi|\sqrt{|\varphi|}}\vartheta^3, \quad \textup{and} \quad  	g(\ddot{V},\dot{x})=\frac{\varepsilon\varepsilon_\lambda d\vartheta}{\alpha\sqrt{|\varphi|}}.
				\label{eq_geo3}
			\end{equation} 
				\begin{enumerate}
					\item If $\lambda$ is timelike or spacelike, then $g(\ddot{V},\dot{x})\neq 0$. Using the fact that $g(\dot{V},\dot{x})$ is a constant, we obtain $g(\ddot{x},\dot{V}) +g(\ddot{V},\dot{x})=0$. Combining the latter identity with equations (\ref{eq_geo3}), we get  $d\neq a+c$, $\frac{\varepsilon_\lambda|\varphi|}{a+c-d}>0$,  and
					\begin{equation}
						\vartheta=\pm\sqrt{\frac{\varepsilon_\lambda|\varphi|}{a+c-d}}.
					\end{equation}
					Additionally, since $g(\dot{x}, \dot{V})=0$, we infer that either $\dot{V}=0$ or $\{V,\dot{V}\}$ forms an orthogonal frame along $x$ and $\dot{x}=\frac{\varepsilon\vartheta}{\sqrt{|\varphi|}}V$. Assuming $\dot{V}=0$, the second equation of (\ref{eq_geo22}) would imply $\vartheta=0$ or $d=0$, which contradicts our assumptions. Therefore, $\dot{V}\neq 0$, $\vartheta\neq 0$, and  $V=\frac{\varepsilon \sqrt{|\varphi|}}{\vartheta}\dot{x}$. Given that $V$ is a unitary vector field, we have   
					\begin{equation*}
						1=\frac{\varepsilon_\lambda|\varphi|}{(a+c)\vartheta^2}=\frac{a+c-d}{a+c},
					\end{equation*}
					which implies $d=0$, resulting in a contradiction. 
					\item If $\lambda$ is lightlike, then $g(\ddot{V},\dot{x})=0$, and equation \eqref{eq_geo3} implies $a+c=d$, leading to $\ddot{x}=0$. By similar reasoning as in (a), we conclude that  $V=\frac{\varepsilon\sqrt{|\varphi|}}{\vartheta}\dot{x}$. However, this leads to 
					$$1=g(V,V)=\frac{|\varphi|}{\vartheta^2}g(\dot{x},\dot{x})=0,$$
					which is impossible.
			\end{enumerate}
		\end{enumerate}
		\cqfd
		
		
		\subsection*{Frenet non null helices on $T_1M^2(\kappa)$}
		
		
		Let $\lambda$ be an arclength-parameterized Frenet non-null curve in $T_1M^2(\kappa)$. Consider the Frenet frame $\lbrace T,W_1,W_2\rbrace$ along $\lambda$. The Serret-Frenet equations can then be written as follows: 
		\begin{equation}
			\left\{\begin{split}
				\tilde{\nabla}_TT=&\varepsilon_1 \kappa_\lambda W_1,\\
				\tilde{\nabla}_TW_1&=-\varepsilon_\lambda \kappa_\lambda T-\varepsilon_2\tau W_2,\\
				\tilde{\nabla}_TW_2&=\varepsilon_1\tau W_1,
			\end{split}\right.
			\label{Frenet_eq}
		\end{equation}
		where  $T=\lambda^\prime$, $\varepsilon_\lambda=\tilde{G}(T,T)$, $\varepsilon_1=\tilde{G}(W_1,W_1)$, and $\varepsilon_2=\tilde{G}(W_2,W_2)$. Here, $\kappa_\lambda$ is the curvature of $\lambda$ and $\tau$ is the torsion of $\lambda$. The quantities $\varepsilon_\lambda$, $\varepsilon_1$ and $\varepsilon_2$ are referred to as \emph{first, second and third causal characters of $\lambda$}, respectively.  
		
		To study helix Frenet non-null curves, we begin with the following theorem, which demonstrates that any helix Frenet curve, with direction the geodesic vector field in $T_1M^2(\kappa)$, satisfies a specific differential equation:
		
		\begin{Theorem}\label{cov_der}
			Let $\lambda:I\subset\mathbb{R}\rightarrow T_1M^2(\kappa)$ be a spacelike or timelike Frenet curve parameterized by its arc length, with first, second, and third causal characters denoted as $\varepsilon_\lambda$, $\varepsilon_1$, and $\varepsilon_2$, respectively. Let $\kappa_\lambda$ and $\tau$ represent its curvature and torsion, and let $\vartheta=\tilde{G}(\lambda',\tilde{\xi})$. 
			\begin{enumerate}
				\item If $\tau=0$, then $\lambda$ is a helix curve with direction along the geodesic vector field in $T_1M^2(\kappa)$ if, and only if, it is a solution of the differential equation
				\begin{equation}\label{diff-eq-helix0}
					\lambda^{(3)} -\frac{\kappa_\lambda'}{\kappa_\lambda}\lambda'' +\varepsilon_1\varepsilon_\lambda\kappa_\lambda^2 \lambda' =0,
				\end{equation}
				or equivalently, the base curve $x$ and the vector field $V$ along $x$ satisfy the system of differential equations
				\label{theorem_sl_cur1}
				\begin{equation}
					\begin{split}
						\nabla_{\dot{x}}\ddot{x}-\frac{\kappa_\lambda'}{\kappa_\lambda}\ddot{x} +\left[\varepsilon_1\varepsilon_\lambda\kappa_\lambda^2 +\frac{d(a+c-d)\vartheta^2}{2\alpha|\varphi|}\right]\dot{x} -\frac{3\varepsilon(a+c-d)\vartheta}{2(a+c)\sqrt{|\varphi|}}\ddot{V}& \\ 
						+\left[-\frac{\varepsilon \kappa_\lambda'(d-(a+c))}{\kappa_\lambda\sqrt{|\varphi|}(a+c)}\vartheta +g(\dot{V},\dot{x})\right]\dot{V} 
						-\left[\frac{\kappa_\lambda'}{\kappa_\lambda}g(\dot{V},\dot{x}) -\frac{3}{2}(g(\ddot{x},\dot{V})+g(\dot{x},\ddot{V}))\right. &\\
						\left.+\frac{\varepsilon( a+c-2d)\vartheta}{2(a+c)\sqrt{|\varphi|}}g(\dot{V},\dot{V})
						-\frac{\varepsilon d^2\vartheta^3}{2\alpha|\varphi|\sqrt{|\varphi|}}+\frac{\varepsilon d\vartheta}{2a\sqrt{|\varphi|}}g(\dot{x},\dot{x})\right]V&=0,
					\end{split}	
					\label{sl_cur100}
				\end{equation} 
				
				\begin{equation}
					\begin{split}
						\frac{\varepsilon  d \kappa_\lambda'\vartheta}{a\kappa_\lambda\sqrt{|\varphi|}}\dot{x}+	\left[\varepsilon_1\varepsilon_\lambda\kappa_\lambda^2 -\frac{(a+c-d)^2\vartheta^2}{2\alpha|\varphi|}+\frac{d\vartheta^2}{a|\varphi|} +g(\dot{V},\dot{V})\right]\dot{V}&\\
						-\frac{\kappa_\lambda'}{\kappa_\lambda}\ddot{V}+\frac{\varepsilon(a+c-3d)}{2a\sqrt{|\varphi|}}\vartheta\ddot{x} +\frac{\varepsilon(a+c-d)\vartheta}{2a\sqrt{|\varphi|}}g(\dot{V},\dot{x})V+\nabla_{\dot{x}}\ddot{V}&=0.
					\end{split}
					\label{sl_cur200}
				\end{equation}   
				\item If $\tau$ doesn't vanish, then $\lambda$ is a helix curve with direction along the geodesic vector field in $T_1M^2(\kappa)$ if, and only if, it is solution of the differential equation 
				\begin{equation}\label{diff-eq-helix}
					\frac{\varepsilon_1\varepsilon_2}{\kappa_\lambda\tau}f_2\lambda^{(3)} -\frac{\varepsilon_1\varepsilon_2\kappa_\lambda'}{\kappa_\lambda^2\tau}f_2\lambda'' -\left[\varepsilon_\lambda\vartheta-\frac{\varepsilon_\lambda\varepsilon_2 }{\tau}\kappa_\lambda f_2\right] \lambda' +\frac{1}{\sqrt{|\varphi|}}V^h=0,
				\end{equation} 
				where $f_2$ is given by the equation 
				\begin{equation}
					f_2=-\frac{\varepsilon_\lambda\vartheta \kappa_\lambda}{\tau} +\frac{\varepsilon_1\varphi}{2\kappa_\lambda\tau\sqrt{|\varphi|}}[g(\ddot{x},\dot{V})-g(\ddot{V},\dot{x})] -\frac{\varepsilon_1\vartheta}{\kappa_\lambda\tau}g(\dot{V},\dot{V})+\frac{\varepsilon_1d\vartheta}{2\alpha \kappa_\lambda\tau}\left(\varepsilon_\lambda-\varepsilon\vartheta^2\right),
					\label{equat1}
				\end{equation}
				or equivalently, the base curve $x$ and the vector field $V$ along $x$ satisfy the system of differential equations
				\begin{equation}
					\begin{split}
						\frac{1}{\kappa_\lambda\tau}f_2\nabla_{\dot{x}}\ddot{x}-\frac{\kappa_\lambda'}{\kappa_\lambda^2\tau}f_2 \ddot{x} +\left[-\varepsilon_1\varepsilon_2\varepsilon_\lambda\vartheta+\frac{1}{\tau}\left(\varepsilon_\lambda \varepsilon_1\kappa_\lambda +\frac{d(a+c-d)\vartheta^2}{2\alpha|\varphi|\kappa_\lambda}\right)f_2\right]\dot{x} & \\ -\frac{3\varepsilon(a+c-d)\vartheta}{2(a+c)\kappa_\lambda\tau\sqrt{|\varphi|}}f_2\ddot{V}
						+\frac{1}{\kappa_\lambda\tau}\left(\frac{\varepsilon \kappa_\lambda'(d-(a+c))}{\kappa_\lambda\sqrt{|\varphi|}(a+c)}\vartheta +g(\dot{V},\dot{x})\right)f_2\dot{V} &\\
						-\left[-\frac{\varepsilon_1\varepsilon_2}{\sqrt{|\varphi|}}+\frac{\kappa_\lambda'}{\kappa_\lambda^2\tau}f_2g(\dot{V},\dot{x}) -\frac{3f_2}{2\kappa_\lambda\tau}(g(\ddot{x},\dot{V})+g(\dot{x},\ddot{V}))\right.\qquad\qquad\qquad &\\
						\left.-\frac{f_2}{\kappa_\lambda\tau}\left(-\frac{\varepsilon( a+c-2d)\vartheta}{2(a+c)\sqrt{|\varphi|}}g(\dot{V},\dot{V})
						+\frac{\varepsilon d^2\vartheta^3}{2\alpha|\varphi|\sqrt{|\varphi|}}-\frac{\varepsilon d\vartheta}{2a\sqrt{|\varphi|}}g(\dot{x},\dot{x})\right)\right]V&=0,
					\end{split}	
					\label{sl_cur1}
				\end{equation} 
				
				\begin{equation}
					\begin{split}
						\frac{\varepsilon  d \kappa_\lambda'\vartheta}{a\kappa_\lambda^2\tau\sqrt{|\varphi|}}f_2\dot{x}-\varepsilon_\lambda \varepsilon_1	\left[\varepsilon_2\vartheta-\frac{\kappa_\lambda }{\tau} f_2+\frac{\varepsilon_\lambda \varepsilon_1}{\kappa_\lambda\tau}f_2\left(\frac{(a+c-d)^2\vartheta^2}{2\alpha|\varphi|}-\frac{d\vartheta^2}{a|\varphi|} -g(\dot{V},\dot{V})\right)\right]\dot{V}&\\
						-\frac{\kappa_\lambda'}{\kappa_\lambda^2\tau}f_2\ddot{V}+\frac{1}{\kappa_\lambda\tau}f_2\left[\frac{\varepsilon(a+c-3d)}{2a\sqrt{|\varphi|}}\vartheta\ddot{x} +\frac{\varepsilon(a+c-d)\vartheta}{2a\sqrt{|\varphi|}}g(\dot{V},\dot{x})V+\nabla_{\dot{x}}\ddot{V}\right]&=0.
					\end{split}
					\label{sl_cur2}
				\end{equation} 
			\end{enumerate}
		\end{Theorem}
		\emph{Proof of Theorem \ref{cov_der}:}
			Deriving the function $\vartheta=\tilde{G}(T,\tilde{\xi})$, we get
			\begin{equation}\label{thet-const}
				\vartheta^\prime=\tilde{G}(\tilde{\nabla}_TT,\tilde{\xi}) +\tilde{G}(\tilde{\nabla}_T\tilde{\xi},T).
			\end{equation}
			On the other hand, since $(\tilde{\eta},\tilde{G},\tilde{\phi},\tilde{\xi})$ is a Sasakian contact pseudo-metric (resp. parasasakian ($\varepsilon$)-paracontact metric) structure, i.e. $h=0$, then we have, using the second identity of\eqref{contact-sasaki},  
			\begin{equation}\label{contact-h=0}
				\tilde{\nabla}_X\tilde{\xi}=-\varepsilon\tilde{\phi} X, \quad \textup{for all}\quad X \in TT_1M^2(\kappa).
			\end{equation}
			Using the first identity of \eqref{Frenet_eq}, \eqref{thet-const}, \eqref{contact-h=0} and \eqref{contact-eta}, we have
			\begin{equation}
				\begin{split}
					G(\tilde{\xi},W_1)= & \frac{\varepsilon_1}{\kappa_\lambda}\tilde{G}(\tilde{\nabla}_TT,\tilde{\xi}) =-\frac{\varepsilon_1}{\kappa_\lambda}\tilde{G}(\tilde{\nabla}_T\tilde{\xi},T) +\frac{\varepsilon_1}{\kappa_\lambda} \vartheta^\prime \\
					= & \frac{\varepsilon\varepsilon_1}{\kappa_\lambda} \tilde{G}(\tilde{\phi}(T),T) +\frac{\varepsilon_1}{\kappa_\lambda} \vartheta^\prime = \frac{\varepsilon\varepsilon_1}{\kappa_\lambda} d\eta(T,T)+\frac{\varepsilon_1}{\kappa_\lambda} \vartheta^\prime= \frac{\varepsilon_1}{\kappa_\lambda} \vartheta^\prime.
				\end{split}
				\label{sl-cur}	
			\end{equation} 
			We deduce that the vector field $\tilde{\xi}$ can be written as 
			\begin{equation}
				\tilde{\xi}=\varepsilon_\lambda\vartheta T+\frac{\vartheta^\prime}{\kappa_\lambda}W_1 +f_2W_2 .
				\label{slant_cur}
			\end{equation}
			\textbf{Case 1:} $\tau=0$. Then \eqref{Frenet_eq} reduces to
			\begin{equation}
				\left\{\begin{split}
					\tilde{\nabla}_TT=&\varepsilon_1 \kappa_\lambda W_1,\\
					\tilde{\nabla}_TW_1&=-\varepsilon_\lambda \kappa_\lambda T.
				\end{split}\right.
			\end{equation}
			We deduce that 
			\begin{equation*}
					\lambda^{(3)} =\varepsilon_1[\kappa_\lambda^\prime W_1 -\varepsilon_\lambda\kappa_\lambda^2T] =\frac{\kappa_\lambda^\prime}{\kappa_\lambda}\lambda^{\prime\prime}-\varepsilon_1\varepsilon_\lambda\kappa_\lambda^2 \lambda^\prime,
			\end{equation*}
			which no other than \eqref{diff-eq-helix0}. To prove that \eqref{diff-eq-helix0} is equivalent to the couple of equations \eqref{sl_cur100} and \eqref{sl_cur200}, we should calculate $\lambda^\prime$, $\lambda^{\prime\prime}$ and $\lambda^{(3)}$. $\lambda^\prime$ is given by \eqref{vel-0}. To calculate $\lambda^{\prime\prime}=\tilde{\nabla}_TT$ and $\lambda^{(3)}=\tilde{\nabla}_T\tilde{\nabla}_TT$, we apply Proposition \ref{cov-der} to get
			\begin{equation}\label{Der2-Lambda}
				\tilde{\nabla}_TT= \left\{\ddot{x} +\frac{d-(a+c)}{a+c} g(\dot{x},V) \dot{V}++g(\dot{x},\dot{V})V\right\}^h +\left\{\ddot{V} -\frac{d}{a}g(\dot{x},V)\dot{x}\right\}^t,
			\end{equation}
			\begin{equation}\label{helix-eq5}
				\begin{split}
					\tilde{\nabla}_T\tilde{\nabla}_TT= & \left\{\nabla_{\dot{x}}\ddot{x} -\frac{d(d-(a+c))}{2\alpha|\varphi|}\vartheta^2\dot{x} +\frac{3\varepsilon(d-(a+c))}{2(a+c)\sqrt{|\varphi|}}\vartheta\ddot{V} +g(\dot{x},\dot{V})\dot{V}\right. \\
					& \left.+\left[\frac{3}{2}\dot{x}(g(\dot{x},\dot{V})) -\frac{\varepsilon d \vartheta}{2a\sqrt{|\varphi|}} g(\dot{x},\dot{x}) +\frac{\varepsilon (2d -(a+c))}{2(a+c)\sqrt{|\varphi|}}\vartheta g(\dot{V},\dot{V}) +\frac{\varepsilon d^2 \vartheta^3}{2\alpha |\varphi|^{\frac{3}{2}} } \right]V\right\}^h \\
					& +\left\{\frac{\varepsilon(a+c-3d)}{2a\sqrt{|\varphi|}}\vartheta\ddot{x} +\nabla_{\dot{x}}\ddot{V} +\frac{\varepsilon(a+c-d)}{2a \sqrt{|\varphi|}}\vartheta g(\dot{x},\dot{V})V\right. \\
					& \left.- \left[g(\ddot{V},V) +\frac{(a+c-d)^2 -2(a+c)d}{2\alpha |\varphi|} \vartheta^2\right]\dot{V}\right\}^t.
				\end{split}
			\end{equation}
			Substituting from \eqref{vel-0}, \eqref{Der2-Lambda} and \eqref{helix-eq5} into the left hand side of \eqref{diff-eq-helix0}, then the horizontal and tangential components are exactly the left hand sides of \eqref{sl_cur100} and \eqref{sl_cur200}, respectively. So \eqref{diff-eq-helix0} is equivalent to the system constituted by \eqref{sl_cur100} and \eqref{sl_cur200}. \\
			\textbf{Case 2:} $\tau$ doesn't vanish. Using the second equation of (\ref{Frenet_eq}), we find that
			$$ W_2 =-\frac{\varepsilon_2}{\tau}[\tilde{\nabla}_TW_1+\varepsilon_\lambda \kappa_\lambda T].$$
			Substuting from the last equation into (\ref{slant_cur}), we obtain 
			\begin{equation}\label{lambda1}
				f_2\tilde{\nabla}_TW_1=-\varepsilon_2\tau \tilde{\xi} +\varepsilon_\lambda[\varepsilon_2\tau \vartheta -\kappa_\lambda f_2]T +\frac{\varepsilon_2 \tau}{\kappa_\lambda}\vartheta^\prime W_1.
			\end{equation}
			On the other hand, deriving the first equation of of (\ref{Frenet_eq}), we get
			$$\tilde{\nabla}_T\tilde{\nabla}_TT=\varepsilon_1[\kappa_\lambda^\prime W_1+\kappa_\lambda\tilde{\nabla}_TW_1] =\frac{\kappa_\lambda^\prime}{\kappa_\lambda}\tilde{\nabla}_TT +\varepsilon_1\kappa_\lambda\tilde{\nabla}_TW_1.$$
			It follows that
			\begin{equation}\label{lambda2}
				\tilde{\nabla}_TW_1=\frac{\varepsilon_1}{\kappa_\lambda}\tilde{\nabla}_T\tilde{\nabla}_TT -\frac{\varepsilon_1\kappa_\lambda^\prime}{\kappa_\lambda^2} \tilde{\nabla}_TT.
			\end{equation}
			Comparing \eqref{lambda1} and \eqref{lambda2}, we get 
			\begin{equation}\label{diff-eq-helix-prime}
				\frac{\varepsilon_1\varepsilon_2}{\kappa_\lambda\tau}f_2\lambda^{(3)} -\frac{\varepsilon_1\varepsilon_2}{\kappa_\lambda^2\tau}[\kappa_\lambda'f_2 +\varepsilon_2\tau\vartheta^\prime]\lambda'' -\left[\varepsilon_\lambda\vartheta-\frac{\varepsilon_\lambda\varepsilon_2 }{\tau}\kappa_\lambda f_2\right] \lambda' +\frac{1}{\sqrt{|\varphi|}}V^h=0.
			\end{equation}
			
			By the generalized fundamental theorem of the local theory of the curves on a manifold, given two real functions $\kappa_\lambda$ and $\tau$, the unique curve, up to a local isometry, with $\kappa_\lambda$ and $\tau$ as curvature and torsion, respectively, is determined by \eqref{diff-eq-helix-prime}. In particular, up to a local isometry, any helix curve is determined by \eqref{diff-eq-helix-prime} with $\vartheta$ is constant, which no other than \eqref{diff-eq-helix}. 
			
			To find the expression of $f_2$ in this case, we begin by differentiating (\ref{sl-cur}), to find
			\begin{equation}\label{helix-eq1}
				\tilde{G}(\nabla_T\tilde{\xi},W_1)+\tilde{G}(\tilde{\xi},\nabla_TW_1)=0.
			\end{equation}
			Using the second identity of \eqref{Frenet_eq} and \eqref{slant_cur}, we obtain 
			\begin{equation}\label{helix-eq2}
				\tilde{G}(\tilde{\xi},\nabla_TW_1)=-\varepsilon_\lambda \kappa_\lambda\tilde{G}(\tilde{\xi},T)+\varepsilon_2\tau\tilde{G}(\tilde{\xi}, W_2)=- \varepsilon_\lambda\vartheta\kappa_\lambda-\tau f_2.
			\end{equation}
			On the other hand, using \eqref{contact-h=0} and \eqref{phi}, we get
			\begin{equation}\label{helix-eq3}
				\tilde{\nabla}_T\tilde{\xi}=-\varepsilon\tilde{\phi}(T) =\frac{\varphi}{2\alpha\sqrt{|\varphi|}}\left[a\dot{V}^h-(a+c)\dot{x}^t\right].
			\end{equation}
			Using equations \eqref{Der2-Lambda}, \eqref{helix-eq3} and the identities $\varepsilon_\lambda=(a+c)g(\dot{x},\dot{x})+\frac{d\vartheta^2}{|\varphi|}+ag(\dot{V},\dot{V})$ and $g(\ddot{V},V)=-g(\dot{V},\dot{V})$, we find 
			\begin{equation}\label{helix-eq4}
				\tilde{G}(\tilde{\nabla}_T\tilde{\xi},W_1)=\frac{\varepsilon_1\varphi}{2\kappa_\lambda\sqrt{|\varphi|}}[g(\ddot{x},\dot{V})-g(\ddot{V},\dot{x})] -\frac{\varepsilon_1\vartheta}{\kappa_\lambda}g(\dot{V},\dot{V})+\frac{\varepsilon_1d\vartheta}{2\alpha \kappa_\lambda}\left(\varepsilon_\lambda-\varepsilon\vartheta^2\right).
			\end{equation}
			Then substituting from \eqref{helix-eq2} and \eqref{helix-eq4} into \eqref{helix-eq1}, we find the expression (\ref{equat1}) of $f_2$. 
			
			Finally, it remains to prove that \eqref{diff-eq-helix} is equivalent to the couple of equations \eqref{sl_cur1} and \eqref{sl_cur2}. Substituting from \eqref{vel-0}, \eqref{Der2-Lambda} and \eqref{helix-eq5} into the left hand side of \eqref{diff-eq-helix}, then the horizontal and tangential components are exactly the left hand sides of \eqref{sl_cur1} and \eqref{sl_cur2}, respectively. So \eqref{diff-eq-helix} is equivalent to the system constituted by \eqref{sl_cur1} and \eqref{sl_cur2}.
		\cqfd
		
		Using Theorem \ref{cov_der}, we shall prove Theorems \ref{th-hor-hel0}-\ref{hel-obl-fren} giving characterizations of horizontal and oblique Frenet helix curves on the unit tangent bundle $T_1M^2(\kappa)$ of a pseudo-Riemannian surface of constant Gaussian curvature $M^2(\kappa)$. 
		
		\emph{Proof of Theorem \ref{th-hor-hel0}:}
		By hypothesis, we have $\dot{V}:=\nabla_{\dot{x}}V=0$, then \eqref{vel-1} becomes
		\begin{equation}
			\varepsilon_\lambda:=\tilde{G}(\lambda',\lambda')=(a+c)g(\dot{x},\dot{x})+\frac{d\vartheta^2}{|\varphi|}.
			\label{hor_hel100}
		\end{equation}
		By Theorem \ref{cov_der}, $\lambda$ is a helix curve directed by the geodesic vector field on $T_1M^2(\kappa)$ if and only if the following system is satisfied 
		\begin{equation}
			\nabla_{\dot{x}}\ddot{x}-\frac{\kappa_\lambda'}{\kappa_\lambda}\ddot{x} +\left[\varepsilon_1\varepsilon_\lambda\kappa_\lambda^2 +\frac{d(a+c-d)\vartheta^2}{2\alpha|\varphi|}\right]\dot{x}  -\frac{\varepsilon d\vartheta}{2a\sqrt{|\varphi|}}\left[g(\dot{x},\dot{x}) -\frac{d\vartheta^2}{(a+c)|\varphi|}\right]V=0,
			\label{hor_hel101}
		\end{equation} 
		
		\begin{equation}
			\vartheta\left[\frac{2d \kappa_\lambda'\vartheta}{\kappa_\lambda}\dot{x}+(a+c-3d)\ddot{x}\right]=0.
			\label{hor_hel102}
		\end{equation}
		\begin{Lemma}\label{Lem10}
			Under the assumptions of Theorem \ref{th-hor-hel0}, if $\lambda$ is a horizontal helix Frenet curve on $T_1M^2(\kappa)$ with zero torsion, directed by $\tilde{\xi}$ on $T_1M^2(\kappa)$, then $\vartheta \neq 0$, $d \neq 0$, $\kappa_\lambda$ is constant and $x$ is a geodesic curve.
		\end{Lemma}
		\emph{Proof of Lemma \ref{Lem10}:} suppose that $\vartheta=0$, then equation (\ref{hor_hel101}) becomes 
		\begin{equation}
			\nabla_{\dot{x}}\ddot{x}-\frac{\kappa_\lambda'}{\kappa_\lambda}\ddot{x} +\varepsilon_1\varepsilon_\lambda\kappa_\lambda^2 \dot{x} =0.
		\end{equation} 	
		Suppose that $\ddot{x}=0$ on some interval, then $\varepsilon_1\varepsilon_\lambda\kappa_\lambda^2 \dot{x} =0$, a contradiction. Then $\ddot{x} \neq 0$ almost everywhere. We deduce that $\{\dot{x},\ddot{x}\}$ is a moving frame along $x$, and hence $V$ can be expressed as $V=\mu \dot{x} +\nu \ddot{x}$. But since $g(\dot{x},V)=0$, $\dot{V}=0$ and $g(\dot{x},\dot{x})$ is constant, then $g(\dot{x},\ddot{x})=g(\ddot{x},V)=0$. Thus $\nu g(\ddot{x},\ddot{x})=0$ and $\mu g(\dot{x},\dot{x})=0$. Since $\mu$ and $\nu$ are not simultaneously zero, then either $g(\dot{x},\dot{x})=0$ or $g(\ddot{x},\ddot{x})=0$. But $g(\dot{x},\dot{x})=\frac{\varepsilon_\lambda}{a+c} \neq 0$, then $g(\ddot{x},\ddot{x})=0$ and $\mu=0$. We deduce that $V=\nu \ddot{x}$. This yields that $1=g(V,V)=\nu^2 g(\ddot{x},\ddot{x})=0$, a contradiction. Then $\vartheta \neq 0$ and  equation \eqref{hor_hel102} becomes
		\begin{equation}
			\frac{2d \kappa_\lambda'\vartheta}{\kappa_\lambda}\dot{x}+(a+c-3d)\ddot{x}=0.
			\label{hor_hel103}
		\end{equation}
		
		Now, making the scalar product of \eqref{hor_hel103} by $\dot{x}$, we obtain $d\kappa_\lambda^\prime=0$. Then either $d=0$ or $\kappa_\lambda^\prime=0$. Suppose now that $d=0$, then equation (\ref{hor_hel103}) becomes $(a+c)\ddot{x}=0$, i.e. $x$ is a geodesic. Hence equation \eqref{hor_hel101} reduces to $\varepsilon_1\varepsilon_\lambda\kappa_\lambda^2 \dot{x} =0$, a contradiction. Then $d \neq 0$ and $\kappa_\lambda^\prime =0$. 
		
		It remains to prove that $x$ is a geodesic curve. Suppose that $x$ is not geodesic, i.e. $\ddot{x} \neq 0$ on an open interval. Then, as before, $\{\dot{x},\ddot{x}\}$ is an orthogonal moving frame along $x$ and  $V$ can be expressed as $V=\mu \dot{x} +\nu \ddot{x}$. Making the scalar product of the last identity by $\ddot{x}$ and using the fact that $g(\dot{x},\ddot{x})=g(V,\ddot{x})=0$, we find that $\nu g(\ddot{x},\ddot{x})=0$. But $g(\ddot{x},\ddot{x}) \neq 0$, since otherwise $g$ would be degenerate. Then $\nu=0$, hence $V=\mu \dot{x}$. Since $g(V,V)=1$ and $g(\dot{x},V)$ are constant, then $\mu$ is a non zero constant. Deriving the expression $V=\mu \dot{x}$, we deduce that $0=\dot{V}=\mu \ddot{x}$, i.e. $\ddot{x}=0$, which is a contradiction. Then $x$ is a geodesic curve.	\cqfd
		
		Suppose that $\lambda$ is a helix curve with zero torsion directed by the geodesic vector field. Then by Lemma \ref{Lem10}, equation \eqref{hor_hel103} is satisfied. On the other hand, equation (\ref{hor_hel101}) has the form  
		\begin{equation}
			A_{01}\dot{x}=A_{02}V,
			\label{hor_hel104}
		\end{equation}
		where 
		\begin{equation}
			A_{01}:=\varepsilon_1\varepsilon_\lambda\kappa_\lambda^2 +\frac{d(a+c-d)\vartheta^2}{2\alpha|\varphi|},
			\label{hor_hel105}
		\end{equation}
		\begin{equation}
			A_{02}:=\frac{\varepsilon d\vartheta}{2\alpha\sqrt{|\varphi|}}\left[\varepsilon_\lambda -\frac{2d\vartheta^2}{|\varphi|}\right].
			\label{hor_hel106}
		\end{equation}
		From \eqref{hor_hel104}, either $A_{01}=A_{02}=0$ or $A_{01}.A_{02} \neq 0$. \\
		\textbf{Case 1:} $A_{01}=A_{02}=0$. $A_{02}=0$ gives $\varepsilon_\lambda=\frac{2d\vartheta^2}{|\varphi|}$. Substituting from the last identity into $A_{01}=0$, we obtain $\kappa_\lambda^2=\frac{\varepsilon_1(d-(a+c))}{4\alpha}$.\\
		\textbf{Case 2:} $A_{01}.A_{02} \neq 0$. Making the scalar product of (\ref{hor_hel104}) by $V$ and by $\ddot{x}$, we obtain 
		\begin{equation}
			\frac{\varepsilon \vartheta}{\sqrt{|\varphi|}}A_{01} =A_{02}, \qquad \frac{1}{a+c}\left(\varepsilon_\lambda -\frac{d\vartheta^2}{|\varphi|}\right)A_{01} =\frac{\varepsilon \vartheta}{\sqrt{|\varphi|}}A_{02}.
			\label{hor_hel107}
		\end{equation}
		Substituting $A_{02}$ from the first equation of \eqref{hor_hel107} into the second equation, we find
		\begin{equation*}
			A_{01}(\varepsilon_\lambda -\varepsilon \vartheta^2)=0.
		\end{equation*}
		Since $A_{01} \neq 0$, then $\vartheta^2 =\varepsilon_\lambda \varepsilon$. We deduce that $\varepsilon_\lambda= \varepsilon$ and $\vartheta^2=1$. Then, it follows then from the first equation of \eqref{hor_hel107}, \eqref{hor_hel105} and \eqref{hor_hel106} that $\kappa_\lambda=0$, which is a contradiction. Then only the first case is valid.
		
		Conversely, it is easy to see that if a non-null horizontal Frenet curve of zero torsion of $T_1M^2(\kappa)$ satisfies the four conditions of the theorem, then it satisfies equations \eqref{hor_hel101} and \eqref{hor_hel102} and is consequently a helix with direction the geodesic vector field.
		\cqfd
		
		\emph{Proof of Theorem \ref{th-hor-hel}:}
		By hypothesis, we have $\dot{V}:=\nabla_{\dot{x}}V=0$, then \eqref{vel-1} becomes
		\begin{equation}
			\varepsilon_\lambda:=\tilde{G}(\lambda',\lambda')=(a+c)g(\dot{x},\dot{x})+\frac{d\vartheta^2}{|\varphi|}.
			\label{hor_hel1}
		\end{equation}
		By Theorem \ref{cov_der}, $\lambda$ is a helix curve directed by the geodesic vector field on $T_1M^2(\kappa)$ if and only if the following system is satisfied 
		\begin{equation}
			\begin{split}
				\frac{1}{\kappa_\lambda\tau}f_2\nabla_{\dot{x}}\ddot{x}-\frac{\kappa_\lambda'}{\kappa_\lambda^2\tau}f_2 \ddot{x}+\left[-\varepsilon_\lambda\varepsilon_1\varepsilon_2\vartheta+\frac{1}{\tau}\left(\varepsilon_\lambda \varepsilon_1\kappa_\lambda +\frac{d(a+c-d)\vartheta^2}{2\alpha|\varphi|\kappa_\lambda}\right)f_2\right]\dot{x}&\\
				+\left[\frac{\varepsilon_1\varepsilon_2}{\sqrt{|\varphi|}}
				+\frac{1}{\kappa_\lambda\tau}\left(\frac{\varepsilon d^2\vartheta^3}{2\alpha|\varphi|\sqrt{|\varphi|}}-\frac{\varepsilon d\vartheta}{2a\sqrt{|\varphi|}}g(\dot{x},\dot{x})\right)f_2\right]V&=0
			\end{split}	
			\label{hor_hel2}
		\end{equation} 
		and
		\begin{equation}
			\left[2d \kappa_\lambda'\dot{x}
			+(a+c-3d)\kappa_\lambda \ddot{x}\right]\vartheta f_2=0,
			\label{hor_hel3}
		\end{equation}
		where \begin{equation}
			f_2=-\varepsilon_\lambda\frac{\vartheta}{\tau}\kappa_\lambda+\frac{\varepsilon_1d\vartheta}{2\alpha \kappa_\lambda\tau}\left(\varepsilon_\lambda-\varepsilon \vartheta^2\right).
			\label{hor_hel4}
		\end{equation}
		\begin{Lemma}\label{Lem1}
			Under the assumptions of Theorem \ref{th-hor-hel}, if $\lambda$ is a horizontal helix Frenet curve on $T_1M^2(\kappa)$, with direction the geodesic vector field on $T_1M^2(\kappa)$, then $\vartheta \neq 0$, $d \neq 0$, $f_2 \neq 0$, i.e. $\varepsilon_2(\varepsilon -\varepsilon_\lambda\vartheta^2) >0$ and $\kappa_\lambda$ is constant.
		\end{Lemma}
		\emph{Proof of Lemma \ref{Lem1}:} suppose that $\vartheta=0$, then $f_2=0$ and equation (\ref{hor_hel2}) becomes $
			\frac{\varepsilon_1\varepsilon_2}{\sqrt{|\varphi|}}V=0,$
		which can not occur, thus $\vartheta\neq 0$. 
		
		On the other hand, equation (\ref{slant_cur}) yields, by making a scalar product operation, 
		\begin{equation}\label{norm}
			f^2_2=\varepsilon_2(\varepsilon-\varepsilon_\lambda\vartheta^2),
		\end{equation}
		hence $\varepsilon_2(\varepsilon-\varepsilon_\lambda\vartheta^2) \geq 0$ and \eqref{hor_hel4} can be written as 
		\begin{equation}\label{hor_hel41}
			f_2=-\varepsilon_\lambda\frac{\vartheta}{\tau}\kappa_\lambda+\frac{\varepsilon_1\varepsilon_2 d\vartheta}{2\alpha \kappa_\lambda\tau}f_2^2.
		\end{equation}
		Suppose that $f_2=0$, then  we deduce from (\ref{hor_hel41}) that $\vartheta\kappa_\lambda=0$, which is a contradiction, since $\kappa_\lambda \neq 0$. Then $f_2 \neq 0$, i.e. $\varepsilon_2(\varepsilon-\varepsilon_\lambda\vartheta^2) > 0$. 
		
		Suppose now that $d=0$, then equation (\ref{hor_hel3}) becomes $\frac{a+c}{2a}\ddot{x}=0$, i.e. $x$ is a geodesic. On the other hand, $f_2=-\frac{\varepsilon_\lambda \kappa_\lambda}{\tau} \vartheta$ and (\ref{hor_hel2}) reduces to 
		\begin{equation}
			V= \sqrt{|a+c|}\vartheta\left[\varepsilon_\lambda+\varepsilon_2\frac{\kappa_\lambda^2}{\tau^2}\right]\dot{x}.
			\label{hor_hel10}
		\end{equation}
		Making the scalar product of (\ref{hor_hel10}) by $\dot{x}$ and using (\ref{hor_hel1}), we obtain $\kappa_\lambda=0$, which is a contradiction, hence $d\neq 0$.
		
		Finally, since $\vartheta \neq 0$ and $f_2 \neq 0$, then equation \eqref{hor_hel3} yields
		\begin{equation*}
			2d \kappa_\lambda'\dot{x} +(a+c-3d)\kappa_\lambda \ddot{x}=0.
			\label{hor_hel31}
		\end{equation*}
			But, since $g(\dot{x},V)$ is constant,  $g(\ddot{x},V)=0$, then
		making the scalar product of the last equation with $V$, we obtain
	 $\kappa_\lambda'=0$, i.e. $\kappa_\lambda$ is constant.	\cqfd
		
		Suppose that $\lambda$ is a helix curve with direction the geodesic vector field. Then by Lemma \ref{Lem1}, equation \eqref{hor_hel3} is equivalent to  
		\begin{equation}
			(a+c-3d)\ddot{x}=0.
			\label{hor_hel14}
		\end{equation}
		On the other hand, equation (\ref{hor_hel2}) has the form  
		\begin{equation}
			A_1\nabla_{\dot{x}}\ddot{x}+A_2\dot{x}+BV=0,
			\label{hor_hel15}
		\end{equation}
		where 
		\begin{equation*}
			A_1:=\frac{1}{\kappa_\lambda \tau}f_2,
			\label{hor_hel16}
		\end{equation*}
		\begin{equation}
			A_2:=-\varepsilon_\lambda\varepsilon_1\varepsilon_2\vartheta+\frac{1}{\tau}\left(\varepsilon_\lambda \varepsilon_1\kappa_\lambda +\frac{d(a+c-d)\vartheta^2}{2\alpha|\varphi|\kappa_\lambda}\right)f_2,
			\label{hor_hel17}
		\end{equation}
		\begin{equation}
			B:=\frac{1}{\sqrt{|\varphi|}}\left[\varepsilon_1\varepsilon_2 -\frac{d\vartheta}{2\alpha \varphi \kappa_\lambda \tau}(\varepsilon_\lambda|\varphi|  -2d \vartheta^2) f_2\right].
			\label{hor_hel18}
		\end{equation}
	Since $ g(\dot{x},\dot{x})$ and $g(V,\dot{x})$ are constants, then making the scalar product of (\ref{hor_hel15}) by $\ddot{x}$, we obtain 
		\begin{equation*}
			A_1g(\nabla_{\dot{x}}\ddot{x},\ddot{x})=0.
		\end{equation*}
		Since $f_2\neq 0$, then $A_1\neq 0$ and hence $g(\nabla_{\dot{x}}\ddot{x},\ddot{x})=0$, and hence $g(\ddot{x},\ddot{x})=\textup{cte}.$
		Let us denote by $K:=\sqrt{|g(\ddot{x},\ddot{x})|}$  the curvature of $x$.
		\begin{itemize}
			\item If $K=0$, then either $x$ is a geodesic or $\ddot{x}$ is a lightlike vector field along $x$. Suppose that  $\ddot{x}$ is a lightlike vector field along $x$, then $\{\dot{x},\ddot{x}\}$ is a local orthogonal frame along $x$, which contradicts the fact that $g$ is non degenerate. Hence $x$ is necessarily a geodesic in this case.
			\item \label{proof} If $K\neq 0$, there exists a unit vector field $W$ along $x$ such that $\ddot{x}=KW$. Since $g(\nabla_{\dot{x}}W,W)=0$, then $\nabla_{\dot{x}}W=\mu\dot{x}$ and then $\nabla_{\dot{x}}\ddot{x}=K\mu \dot{x}$. It follows that 
			\begin{equation*}
					K\mu g(\dot{x},\dot{x})=g(\nabla_{\dot{x}}\ddot{x},\dot{x})=-K^2\frac{
						g(\ddot{x},\ddot{x})}{|	g(\ddot{x},\ddot{x})|}.
			\end{equation*}
			In this case, $g(\dot{x},\dot{x})\neq 0$ and, by virtue of \eqref{hor_hel1}, we deduce that 
			$$\mu=-K\varepsilon'(a+c)\left(\frac{|\varphi|}{\varepsilon_\lambda|\varphi|-d\vartheta^2}\right),$$
			where $\varepsilon'=\frac{
				g(\ddot{x},\ddot{x})}{|	g(\ddot{x},\ddot{x})|}$. 
			It follows that  
			\begin{equation*}
				\nabla_{\dot{x}}\ddot{x}=-\varepsilon'K^2\left(\frac{(a+c)|\varphi|}{\varepsilon_\lambda|\varphi|-d\vartheta^2}\right)\dot{x},
			\end{equation*}
			in other words, $x$ is a pseudo-Riemannian circle in $M^2(\kappa)$.
		\end{itemize}
		In both cases, i.e. $x$ is a geodesic ($K=0$) or $x$ is a pseudo-Riemannian circle ($K\neq 0$), we have 
		\begin{equation}
			\nabla_{\dot{x}}\ddot{x}=K^\prime\dot{x},
			\label{hor_hel21}
		\end{equation}
		where
		\begin{equation}\label{K-prime}
			K^\prime=\left\{
			\begin{array}{ll}
				0, & \textup{if } K=0, \\
				-\varepsilon'K^2\left(\frac{(a+c)|\varphi|}{\varepsilon_\lambda|\varphi|-d\vartheta^2}\right), & \textup{if } K \neq 0.
			\end{array}
			\right.
		\end{equation}
		Substituting from (\ref{hor_hel21}) into (\ref{hor_hel15}), we obtain 
		\begin{equation}
			A\dot{x}+BV=0,
			\label{hor_hel22}
		\end{equation}
		where $A=K^\prime +A_2.$
		Hence the two possible cases are $A=B=0$ or $AB\neq 0$.
		\begin{Lemma}\label{Lem2}
			Under the assumptions of Lemma \ref{Lem1}, we have $A=B=0$.
		\end{Lemma}
		\emph{Proof of Lemma \ref{Lem2}:} Suppose that $A\neq 0$ and $B\neq 0$, then we have, by virtue of \eqref{hor_hel22}, $V=-\frac{A}{B}\dot{x}.$
		Since $g(V,V)=1$ and $g(\dot{x},\dot{x})=\frac{1}{a+c}\left(\varepsilon_\lambda-\frac{d\vartheta^2}{|\varphi|}\right)$, 
		then
		\begin{equation}
			\frac{A}{B}=\pm\sqrt{\frac{(a+c)|\varphi|}{\varepsilon_\lambda|\varphi|-d\vartheta^2}},
			\label{hor_hel25}
		\end{equation}
		and consequently
		\begin{equation}
			V=\pm\sqrt{\frac{(a+c)|\varphi|}{\varepsilon_\lambda|\varphi|-d\vartheta^2}}\dot{x}.
			\label{hor_hel26}
		\end{equation}
		Using the fact that $g(\dot{x},V)=\frac{\varepsilon}{\sqrt{|\varphi|}}\vartheta$, we obtain from (\ref{hor_hel26})
			$\vartheta^2=\frac{\varepsilon_\lambda|\varphi|}{\varphi}=\varepsilon_\lambda\varepsilon.$
		Hence $\varepsilon_\lambda\varepsilon=1$, i.e. $\tilde{\xi}$ and $\lambda$ have the same causal character and $\vartheta^2=1$. We deduce that 
		$\varepsilon_\lambda|\varphi|= \varepsilon|\varphi|=\varphi$ and, substituting into \eqref{hor_hel25}, we find $\frac{A}{B}=\pm\sqrt{|\varphi|}.$
		On the other hand, \eqref{hor_hel4} becomes $f_2=-\varepsilon_\lambda\vartheta\frac{\kappa_\lambda}{\tau}$ and, therefore equation \eqref{norm} yields  $\kappa_\lambda=0$, which is a contradiction. \cqfd
		Since, by Lemma \ref{Lem2}, $A=B=0$, then equation (\ref{hor_hel22}) is automatically satisfied and, by virtue of \eqref{hor_hel18}, $B=0$ implies that 
		\begin{equation}\label{hor_hel231}
			f_2=\frac{2\varepsilon_1\varepsilon_2\alpha\varphi \kappa_\lambda\tau}{d(\varepsilon_\lambda|\varphi| -2d\vartheta^2)\vartheta}.
		\end{equation}
		Comparing \eqref{hor_hel231} with \eqref{hor_hel4}, we find that the function $\tau$ is a constant satisfying
		\begin{equation}\label{hor_torsion}
			\tau^2=\frac{\varepsilon_2d\vartheta^2}{2\alpha\varphi}(\varepsilon_\lambda|\varphi| 
				-2d\vartheta^2)\left[\frac{d(\varepsilon_\lambda -\varepsilon\vartheta^2)}{2\alpha \kappa_\lambda^2}-\varepsilon_\lambda\varepsilon_1 \right],
		\end{equation} 
		on the other hand \eqref{norm} implies 
		\begin{equation} \label{hor_torsion2}
			\tau^2=\frac{\varepsilon_2d^2\vartheta^2(\varepsilon-\varepsilon_\lambda\vartheta^2)}{4\alpha^2\varphi^2 
				\kappa_\lambda^2}(\varepsilon_\lambda|\varphi| -2d\vartheta^2)^2.
		\end{equation}
		We deduce, in particular, that 
			$\varepsilon_2(\varepsilon-\varepsilon_\lambda\vartheta^2)>0.$
		Comparing equations \eqref{hor_torsion} and \eqref{hor_torsion2}, we get 
		\begin{equation}\label{hor_torsion21}
			\kappa_\lambda^2=\frac{\varepsilon_\lambda\varepsilon_1 d^2}{\alpha\varphi}(\varepsilon-\varepsilon_\lambda\vartheta^2)\vartheta^2.
		\end{equation}
		
		Substituting from the last equation into \eqref{hor_torsion2}, we find $
			\tau^2=\frac{(\varepsilon_\lambda|\varphi|-2d\vartheta^2)^2}{4\alpha\varphi}.$
		It follows, in particular, that $\alpha\varphi>0$. Since $4\alpha=\pm|\varphi|$ then $4\alpha=\varphi.$
		
		Then substuting from \eqref{hor_hel231} and \eqref{hor_torsion21} into \eqref{hor_hel17}, we get $A_2=0$ and since $A=0$ then $K^\prime=0,$
		hence $x$ is a geodesic.\\
		Conversely, it is easy to see that if a non-null horizontal Frenet curve of $T_1M^2(\kappa)$ satisfies the four conditions of the theorem, then it satisfies equations \eqref{hor_hel2} and \eqref{hor_hel3} and is consequently a helix with axis the geodesic vector field.
		\cqfd

		\emph{Proof of Theorem \ref{hel-obl-fren0}:}
		The scalar product of equation (\ref{sl_cur200}) by $V$ gives, by virtue of $g(\dot{V},V)=0$, 
		\begin{equation}\label{hel_cur39}
			\frac{\varepsilon  d \kappa_\lambda'\vartheta^2}{a\kappa_\lambda|\varphi|}
			-\frac{\kappa_\lambda'}{\kappa_\lambda}g(\ddot{V},V)+\frac{\varepsilon(a+c-3d)}{2a\sqrt{|\varphi|}}\vartheta g(\ddot{x},V) +\frac{\varepsilon(a+c-d)\vartheta}{2a\sqrt{|\varphi|}}g(\dot{V},\dot{x})+g(\nabla_{\dot{x}}\ddot{V},V)=0.
		\end{equation}
		By \eqref{vel-1}, we have
		\begin{equation}\label{hel_cur03}
			g(\dot{V},\dot{V})=\frac{1}{a}\left[\varepsilon_\lambda-(a+c)\sigma-\frac{d\vartheta^2}{|\varphi|}\right].
		\end{equation}
		Since $g(\dot{x},\dot{x})=\sigma$ is constant, then $g(\dot{V},\dot{V})$ is constant. Differentiating twice $g(\dot{V},V)=0$, we have
		\begin{equation}\label{prod-x-V}
			g(\ddot{V},V)=-g(\dot{V},\dot{V}), \qquad g(\nabla_{\dot{x}}\ddot{V},V)=-\frac32 \frac{d}{dt}(g(\dot{V},\dot{V}))=0.
		\end{equation}
		On the other hand, differentiating twice the constant function $g(\dot{x},V)$, we obtain
		\begin{equation}\label{prod-V-V}
			g(\ddot{x},V)=-g(\dot{x},\dot{V}), \qquad g(\nabla_{\dot{x}}\ddot{x},V)=-2g(\ddot{x},\dot{V})-g(\dot{x},\ddot{V}).
		\end{equation}
		Substituting from \eqref{prod-V-V} and \eqref{prod-x-V} into \eqref{hel_cur39}, we get
		\begin{equation}
			\frac{\kappa_\lambda'}{\kappa_\lambda}\left[g(\dot{V},\dot{V}) +\frac{d \vartheta^2}{a|\varphi|}\right]	
			+\frac{\varepsilon d\vartheta}{a\sqrt{|\varphi|}}g(\dot{V},\dot{x})=0.
			\label{hel_cur30}
		\end{equation}
		We note that $\{V,\dot{V}\}$ is an orthogonal frame along the curve $x$ and $g(\dot{V},\dot{V})\neq 0$, otherwise $g$ will be degenerate.	
		On the other hand, if we write $\dot{x}=\mu V+\nu \dot{V}$, then $\mu=g(\dot{x},V)$ is constant and $\sigma=g(\dot{x},\dot{x})=\mu^2 +\nu^2g(\dot{V},\dot{V})$ is also constant by hypothesis. Then $\nu^2$ is constant, since $g(\dot{V},\dot{V})\neq 0$ is a constant. It follows that $g(\dot{x},\dot{V})=\nu g(\dot{V},\dot{V})$ is constant. More precisely, we have
		\begin{equation}\label{cte-20}
			g(\dot{x},\dot{V})^2=\left(\sigma -\frac{\vartheta^2}{|\varphi|}\right)g(\dot{V},\dot{V}).
		\end{equation}
		Furthermore, $g(\dot{V},\dot{V})$ is constant yields $g(\ddot{V},\dot{V})=0$ and 
		\begin{equation}\label{hel_cur1410}
			\ddot{V}=-g(\dot{V},\dot{V})V, \quad \textup{and} \quad
			\nabla_{\dot{x}}\ddot{V}=-g(\dot{V},\dot{V})\dot{V}.
		\end{equation}
		We also have $g(\ddot{x},\dot{V})$ and $g(\ddot{V},\dot{x})$ are constant. Indeed, making the scalar product of \eqref{hel_cur1410} by $\dot{x}$, we get $g(\ddot{V},\dot{x})=-g(\dot{V},\dot{V})g(\dot{x},V)$, which is constant. On the other hand, since $g(\dot{x},\dot{V})$ is constant then $g(\ddot{x},\dot{V})=-g(\ddot{V},\dot{x})$ and, therefore, $g(\ddot{x},\dot{V})$ is constant. 
		
		 Suppose that $\vartheta=0$. Then we have, by virtue of \eqref{hel_cur30}, $\kappa_\lambda^\prime=0$. Substituting into \eqref{sl_cur200}, we obtain $\varepsilon_1\varepsilon_\lambda\kappa_\lambda^2 \dot{V}=0$, which is a contradiction. Then $\vartheta\neq 0$.
		
		The scalar product of equation (\ref{sl_cur200}) by $\dot{V}$ gives, by virtue of \eqref{hel_cur1410}, 
		\begin{equation}\label{hel_cur391}
			\frac{ d \kappa_\lambda'\vartheta^2}{a\kappa_\lambda|\varphi|}
			+\left[\varepsilon_1\varepsilon_\lambda\kappa_\lambda^2+\frac{d(a+c-d)}{2\alpha |\varphi|}\vartheta^2\right] g(\dot{V},\dot{V})=0.
		\end{equation}
		It follows from the preceding equation that $d \neq 0$. Hence \eqref{hel_cur30} gives
		\begin{equation}
			g(\dot{V},\dot{x})=-\frac{\varepsilon a\sqrt{|\varphi|}\kappa_\lambda'}{d\vartheta\kappa_\lambda}\left[g(\dot{V},\dot{V}) +\frac{d \vartheta^2}{a|\varphi|}\right].
			\label{hel_cur301}
		\end{equation}
		In particular, we have
		\begin{equation}\label{cte-lambda'}
			\frac{\kappa_\lambda^\prime}{\kappa_\lambda} \textup{ is constant}.
		\end{equation}
		
		Using the fact that $\{V,\dot{V}\}$ is an orthogonal moving frame along $x$, \eqref{eq4} and \eqref{hel_cur301} yield
		\begin{equation}\label{expr-x-dot}
			\dot{x}=A_0 V+B_0\dot{V}, 
		\end{equation}
		where
		\begin{equation}\label{A0-B0}
			A_0=\frac{\varepsilon \vartheta}{\sqrt{|\varphi|}}, \qquad B_0=-\frac{\varepsilon a\sqrt{|\varphi|}\kappa_\lambda'}{d\vartheta\kappa_\lambda}\left[1 +\frac{d \vartheta^2}{a|\varphi|g(\dot{V},\dot{V})}\right].
		\end{equation}
		Since $A_0$ and $B_0$ are constant, then deriving twice \eqref{expr-x-dot} and using \eqref{hel_cur1410}, we get
		\begin{equation}\label{expr-x-ddot}
			\ddot{x}=-B_0g(\dot{V},\dot{V}) V+A_0\dot{V}, \qquad \nabla_{\dot{x}}\ddot{x}=-g(\dot{V},\dot{V})(A_0\dot{V} +B_0 \dot{V}).
		\end{equation}
		Substituting from \eqref{expr-x-dot}, \eqref{expr-x-ddot} and \eqref{hel_cur1410}  into \eqref{sl_cur200} and taking the $\dot{V}$-component, we obtain
		\begin{equation}\label{hel_cur1411}
			\frac{\varepsilon d\vartheta B_0 \kappa_\lambda^\prime}{a\sqrt{|\varphi|}\kappa_\lambda} +\varepsilon_1\varepsilon_\lambda \kappa_\lambda^2 +\frac{d(a+c-d)}{2\alpha|\varphi|}\vartheta^2=0.
		\end{equation}
		By virtue of \eqref{cte-lambda'}, we deduce from the last equation that $\kappa_\lambda$ is constant, i.e. $\kappa_\lambda^\prime=0$. Then from the second expression of \eqref{A0-B0}, we have $B_0=0$ and \eqref{expr-x-dot} reduces to
		\begin{equation}\label{expr-x-dot2}
			\dot{x}=A_0 V=\frac{\varepsilon \vartheta}{\sqrt{|\varphi|}}V. 
		\end{equation}
		It follows that 
		\begin{equation}\label{expr-x-dot3}
			\sigma=g(\dot{x},\dot{x})=\frac{\vartheta^2}{|\varphi|} \quad \textup{and} \quad \dot{x}=\pm \sqrt{\sigma}V. 
		\end{equation}
		Since $\kappa_\lambda$ is constant, then \eqref{hel_cur1411} yields
		\begin{equation}\label{kappa0}
			\kappa_\lambda^2= -\frac{\varepsilon_1 \varepsilon_\lambda d(a+c-d)}{2\alpha|\varphi|}.
		\end{equation}
		On the other hand, by \eqref{expr-x-dot3}, we have $\ddot{x}=\pm \sqrt{\sigma}\dot{V}$, hence $g(\ddot{x},\ddot{x})=\sigma g(\dot{V},\dot{V})$ is a non zero constant. We deduce that $x$ is a pseudo-Riemannian circle.
		
		Conversely, if an oblique Frenet curve $\lambda$ with no torsion satisfies the four conditions of the theorem yield \eqref{sl_cur100} and \eqref{sl_cur200}, i.e. $\lambda$ is a helix curve with direction the geodesic vector field.
		\cqfd

		\emph{Proof of Theorem \ref{hel-obl-fren}:}
		Since $g(\dot{x},V)$ and $g(V,V)$ are constant, then $g(\ddot{x},V)=-g(\dot{x},\dot{V})$ and $g(\ddot{V},V)=-g(\dot{V},\dot{V})$.
		Then the scalar products of equation (\ref{sl_cur1}) by $\dot{V}$ and $V$, respectively, give 
		\begin{equation}
			\begin{split}				
				\frac{\varepsilon_1\varepsilon_2}{\kappa_\lambda\tau}f_2g(\nabla_{\dot{x}}\ddot{x},\dot{V}) -\frac{\varepsilon_1\varepsilon_2\kappa_\lambda'}{\kappa_\lambda^2\tau}f_2 g(\ddot{x},\dot{V})
				-\frac{3\varepsilon\varepsilon_1\varepsilon_2(a+c-d)\vartheta}{2(a+c)\kappa_\lambda\tau\sqrt{|\varphi|}}f_2g(\ddot{V},\dot{V})&\\
				+\left[-\varepsilon_\lambda\vartheta+\frac{1}{\tau}\left(\varepsilon_\lambda \varepsilon_2\kappa_\lambda +\frac{\varepsilon_1\varepsilon_2d(a+c-d)\vartheta^2}{2\alpha|\varphi| \kappa_\lambda}\right)f_2\right]g(\dot{x},\dot{V})\\
				+\frac{\varepsilon_1\varepsilon_2}{\kappa_\lambda\tau}\left(\frac{\varepsilon \kappa_\lambda'(d-(a+c))}{\kappa_\lambda\sqrt{|\varphi|}(a+c)}\vartheta -g(V,\ddot{x})\right)f_2g(\dot{V},\dot{V}) &=0
			\end{split}	
			\label{cur_hel1}
		\end{equation} 
		and
		\begin{equation}
			\begin{split}				\frac{\varepsilon_1\varepsilon_2}{\kappa_\lambda\tau}f_2g(\nabla_{\dot{x}}\ddot{x},V)-\frac{1}{\sqrt{|\varphi|}}(\varepsilon\varepsilon_\lambda\vartheta^2-1) +\frac{\varepsilon \vartheta}{\sqrt{|\varphi|}\tau}\left(\varepsilon_\lambda \varepsilon_2\kappa_\lambda +\frac{\varepsilon_1\varepsilon_2d(a+c-d)\vartheta^2}{2\alpha|\varphi|\kappa_\lambda}\right)f_2&\\
				+\frac{\varepsilon_1\varepsilon_2}{\kappa_\lambda\tau}\left(\frac{3}{2}\dot{x}(g(\dot{x},\dot{V}))+\frac{\varepsilon( 2(a+c)-d)\vartheta}{2(a+c)\sqrt{|\varphi|}}g(\dot{V},\dot{V})
				+\frac{\varepsilon d^2\vartheta^3}{2\alpha|\varphi|\sqrt{|\varphi|}}-\frac{\varepsilon d\vartheta}{2a\sqrt{|\varphi|}}\sigma\right)f_2&=0.
			\end{split}	
			\label{cur_hel2}
		\end{equation}
		
		\begin{Lemma}\label{Lem3}
			Under the assumptions of Theorem \ref{hel-obl-fren}, if $\lambda$ is an oblique helix Frenet curve on $T_1M^2(\kappa)$, with direction $\tilde{\xi}$ on $T_1M^2(\kappa)$, then $\vartheta \neq 0$, $d \neq 0$ and $f_2 \neq 0$, i.e. $\varepsilon_2(\varepsilon -\varepsilon_\lambda\vartheta^2) >0$.
		\end{Lemma} 
		\emph{Proof of Lemma \ref{Lem3}:} Suppose that $f_2=0$, then equation (\ref{cur_hel1}) becomes
		\begin{equation*}
			\vartheta g(\dot{x},\dot{V})=0,
		\end{equation*}
		but equation (\ref{sl_cur1}) implies  $\vartheta\neq 0$, then $g(\dot{x},\dot{V})=0$. Making use of equation (\ref{cur_hel2}), we find
       $\vartheta^2=\varepsilon\varepsilon_\lambda=1.$
		Moreover, equation (\ref{sl_cur2}) becomes $\dot{V}=0$, which contradicts the fact that $\lambda$ is oblique, hence $f_2\neq0$.
		
		Let us prove that $\vartheta \neq 0$. Since $f_2 \neq 0$, then the scalar product of equation (\ref{sl_cur2}) by $V$ gives 
		\begin{equation}
			\frac{ d \kappa_\lambda'\vartheta^2}{a\kappa_\lambda|\varphi|}
			+\frac{\kappa_\lambda'}{\kappa_\lambda}g(\dot{V},\dot{V})+\frac{\varepsilon 
				d\vartheta}{a\sqrt{|\varphi|}}g(\dot{V},\dot{x})+g(\nabla_{\dot{x}}\ddot{V},V)=0.
			\label{hel_cur3}
		\end{equation}
		We note that $\{V,\dot{V}\}$ is an orthogonal frame along the curve $x$ and $g(\dot{V},\dot{V})\neq 0$, otherwise $g$ will be degenerate.
		Differentiating twice the expression $g(\dot{V},V)=0$, we obtain $g(\nabla_{\dot{x}}\ddot{V},V)=-3g(\ddot{V},\dot{V})$.	Then \eqref{vel-1} implies
		\begin{equation}
			g(\dot{V},\dot{V})=\frac{1}{a}\left[\varepsilon_\lambda-(a+c)\sigma-\frac{d\vartheta^2}{|\varphi|}\right],
		\end{equation}
		in particular $g(\dot{V},\dot{V})$ is constant, thus $g(\ddot{V},\dot{V})=0$ and $g(\nabla_{\dot{x}}\ddot{V},V)=0$. Hence (\ref{hel_cur3}) becomes 
		\begin{equation}
			\frac{ d \kappa_\lambda'\vartheta^2}{a\kappa_\lambda|\varphi|}
			+\frac{\kappa_\lambda'}{\kappa_\lambda}g(\dot{V},\dot{V})+\frac{\varepsilon d\vartheta}{a\sqrt{|\varphi|}}g(\dot{V},\dot{x})=0.
			\label{hel_cur33}
		\end{equation}	
		On the other hand, if we write $\dot{x}=\mu V+\nu \dot{V}$, then $\mu=g(\dot{x},V)$ is constant and $\sigma=g(\dot{x},\dot{x})=\mu^2 +\nu^2g(\dot{V},\dot{V})$ is also constant by hypothesis. Then $\nu^2$ is constant, since $g(\dot{V},\dot{V})\neq 0$ is a constant. It follows that $g(\dot{x},\dot{V})=\nu g(\dot{V},\dot{V})$ is constant. More precisely, we have
		\begin{equation}\label{cte-2}
			g(\dot{x},\dot{V})^2=g(\dot{V},\dot{V})\left(\sigma -\frac{\vartheta^2}{|\varphi|}\right).
		\end{equation}
		Furthermore, $g(\dot{V},\dot{V})$ is constant yields $g(\ddot{V},\dot{V})=0$ and 
		\begin{equation}\label{hel_cur141}
			\ddot{V}=-g(\dot{V},\dot{V})V,
		\end{equation}
		hence 
		\begin{equation}
			\nabla_{\dot{x}}\ddot{V}=-g(\dot{V},\dot{V})\dot{V}.
			\label{equat2}
		\end{equation}
		We also have $g(\ddot{x},\dot{V})$ and $g(\ddot{V},\dot{x})$ are constant. Indeed, making the scalar product of \eqref{hel_cur141} by $\dot{x}$, we get $g(\ddot{V},\dot{x})=-g(\dot{V},\dot{V})g(\dot{x},V)$, which is constant. On the other hand, since $g(\dot{x},\dot{V})$ is constant then $g(\ddot{x},\dot{V})=-g(\ddot{V},\dot{x})$ and, therefore, $g(\ddot{x},\dot{V})$ is constant. Finally, from the preceding, \eqref{equat1} becomes
		\begin{equation}\label{equat11}
			f_2=-\varepsilon_\lambda\vartheta\frac{\kappa_\lambda}{\tau} +\frac{\varepsilon_1d \vartheta}{2\alpha\kappa_\lambda\tau}(\varepsilon_\lambda -\varepsilon \vartheta^2).
		\end{equation}
		Since $f_2 \neq 0$, then $\vartheta\neq 0$.      
		
		Now, we prove that $d \neq 0$. Suppose that  $d=0$, then equation (\ref{hel_cur33}) yields $\kappa_\lambda'=0$ and equations (\ref{sl_cur1}) and (\ref{sl_cur2}) become, 
		respectively,
		\begin{equation}
			\begin{split}
				\frac{1}{\kappa_\lambda\tau}f_2\nabla_{\dot{x}}\ddot{x}+\varepsilon_\lambda\varepsilon_1\left[-\varepsilon_2\vartheta +\frac{\kappa_\lambda}{\tau}f_2\right]\dot{x} -\frac{3\varepsilon\vartheta}{2\kappa_\lambda\tau\sqrt{|\varphi|}}f_2\ddot{V} &  \\
				-\frac{1}{\kappa_\lambda\tau}g(V,\ddot{x})f_2\dot{V} 
				+\left[\frac{\varepsilon_1\varepsilon_2}{\sqrt{|\varphi|}}-\frac{\varepsilon\vartheta}{2\kappa_\lambda\tau\sqrt{|\varphi|}}g(\dot{V},\dot{V})
				f_2\right]V  & =0
			\end{split}
			\label{sl_cur11}
		\end{equation} 
		and
		\begin{equation}
			\begin{split}
				\left[-\varepsilon_\lambda\varepsilon_1\varepsilon_2\vartheta+\frac{\varepsilon_\lambda \varepsilon_1\kappa_\lambda }{\tau} f_2-\frac{1}{\kappa_\lambda\tau}f_2\left(\frac{(a+c)}{2a|\varphi|}\vartheta^2-g(\dot{V},\dot{V})\right)\right]\dot{V} &  \\
				+\frac{\varepsilon(a+c)\vartheta}{2a\sqrt{|\varphi|}\kappa_\lambda\tau}f_2\left[\ddot{x}+g(\dot{V},\dot{x})V\right] +\frac{1}{\kappa_\lambda\tau}f_2\nabla_{\dot{x}}\ddot{V} & =0.
			\end{split}
			\label{sl_cur21}
		\end{equation}
		The scalar product of equation (\ref{sl_cur11}) by $\dot{V}$ gives 
		\begin{equation}			
			\frac{1}{\kappa_\lambda\tau}f_2g(\nabla_{\dot{x}}\ddot{x},\dot{V})+\varepsilon_\lambda\varepsilon_1\left[-\varepsilon_2\vartheta +\frac{\kappa_\lambda}{\tau}f_2\right]g(\dot{x},\dot{V})-\frac{1}{\kappa_\lambda\tau}g(V,\ddot{x})f_2g(\dot{V},\dot{V})=0.
			\label{hel_cu1}
		\end{equation} 
		Differentiating $g(\dot{x},\dot{V})=\textup{cte}$ and using \eqref{equat2}, we obtain
		\begin{equation}\label{hel_cu11}
			g(\ddot{x},\dot{V})=-g(\dot{x},\ddot{V})=-g(\dot{V},\dot{V})g(\dot{x},V).
		\end{equation}
		
		On the other hand, making a second differentiation of $g(\dot{x},\dot{V})=\textup{cte}$ and using \eqref{hel_cur141}, \eqref{equat2} and the identity $g(\ddot{x},V)=-g(\dot{x},\dot{V})$, we obtain 
		\begin{equation*}
				g(\nabla_{\dot{x}}\ddot{x},\dot{V})=-2g(\ddot{x},\ddot{V})-g(\dot{x},\nabla_{\dot{x}}\ddot{V})
				=g(\ddot{x},V)g(\dot{V},\dot{V}).
			\label{hel_cu2}
		\end{equation*}
		Substituting from the last equation into (\ref{hel_cu1}), we deduce that 
		\begin{equation}\label{hel_cu21}
			\left[-\varepsilon_2\vartheta +\frac{\kappa_\lambda}{\tau}f_2\right]g(\dot{x},\dot{V})=0.
		\end{equation}
		Then either $g(\dot{x},\dot{V})=0$ or $-\varepsilon_2\vartheta +\frac{\kappa_\lambda}{\tau}f_2=0$. We claim that $g(\dot{x},\dot{V}) \neq 0$. Indeed, suppose that $g(\dot{x},\dot{V})=0$. Then, we have $\dot{x}=\frac{\varepsilon\vartheta}{\sqrt{|\varphi|}}V$. Hence, using \eqref{hel_cur141}, equations \eqref{sl_cur11} and \eqref{sl_cur21} become, respectively,
		\begin{equation}
			\left[-\vartheta+\varepsilon_2\frac{\kappa_\lambda}{\tau}f_2\right]+\frac{\varepsilon_\lambda}{\sqrt{|\varphi|}}=0,\quad -\vartheta+\varepsilon_2\frac{\kappa_\lambda}{\tau}f_2=0,
		\end{equation}
		thus $\varepsilon_\lambda=0$, which is a contradiction. We deduce that $g(\dot{x},\dot{V}) \neq 0$ 
		and hence $$-\varepsilon_2\vartheta +\frac{\kappa_\lambda}{\tau}f_2=0.$$ 
		By virtue of \eqref{equat11},
			$f_2=-\varepsilon_\lambda\vartheta\frac{\kappa_\lambda}{\tau},$
	hence $\tau$ is constant satisfying $
			\tau^2=-\varepsilon_\lambda\varepsilon_2 \kappa^2_\lambda.$
		Then using $f^2_2=\varepsilon_2(\varepsilon-\varepsilon_\lambda\vartheta^2)$, we find $\vartheta^2\kappa_\lambda^2=\varepsilon_2(\varepsilon-\varepsilon_\lambda\vartheta^2)\tau^2=\varepsilon\varepsilon_2+\vartheta^2\kappa_\lambda^2$, hence $\varepsilon\varepsilon_2=0$, which is a contradiction, hence $d\neq 0$.\cqfd
		
		Now, equation (\ref{hel_cur33}) implies  that 
		\begin{equation}
			g(\dot{V},\dot{x})=-\frac{\varepsilon a\sqrt{|\varphi|}\kappa_\lambda'}{d\vartheta 
				\kappa_\lambda}\left[\frac{d\vartheta^2}{a|\varphi|}+g(\dot{V},\dot{V})\right]
			\label{hel_cur4}
		\end{equation}
		Hence  $\frac{\kappa_\lambda'}{\kappa_\lambda}$ is constant, since $g(\dot{V},\dot{x})$ is constant, and  
		\begin{equation}
			\dot{x}=\frac{\varepsilon\vartheta}{\sqrt{|\varphi|}}V-\frac{\varepsilon a \sqrt{|\varphi|}\kappa_\lambda'}{d\vartheta \kappa_\lambda g(\dot{V},\dot{V})}\left[\frac{d\vartheta^2}{a|\varphi|}+g(\dot{V},\dot{V})\right]\dot{V}.
			\label{hel_cur5}
		\end{equation}
		Now, from \eqref{eq4} and \eqref{cte-2}, we have 
		\begin{equation*}
			\dot{x}=\frac{\varepsilon\vartheta}{\sqrt{|\varphi|}}V\pm\sqrt{\frac{1}{g(\dot{V},\dot{V})}\left(\sigma-\frac{\vartheta^2}{|\varphi|}\right)}\dot{V},
		\end{equation*}
		then  
		
		\begin{equation}
			g(\dot{x},\dot{V})=\pm\sqrt{g(\dot{V},\dot{V})\left(\sigma-\frac{\vartheta^2}{|\varphi|}\right)}.
			\label{hel_cu3}
		\end{equation}
		Taking into account \eqref{hel_cu3}, \eqref{hel_cur4} yields
		\begin{equation} 
			\frac{\kappa_\lambda'}{\kappa_\lambda}=\pm\frac{d\vartheta|\varphi|}{d\vartheta^2+a|\varphi|g(\dot{V},\dot{V})}\sqrt{\left(\frac{\sigma}{|\varphi|} 
				-\frac{\vartheta^2}{\varphi^2}\right)g(\dot{V},\dot{V})}
			\label{hel_cur6}
		\end{equation}
		On the other hand, we have $g(\ddot{V},\dot{V})=0$ and $g(\ddot{V},V)=-g(\dot{V},\dot{V})$, thus 
		\begin{equation}
			\ddot{V}=-g(\dot{V},\dot{V})V
			\label{hel_cur7}
		\end{equation}
		and hence
		\begin{equation}
			\nabla_{\dot{x}}\ddot{V}=-g(\dot{V},\dot{V})\dot{V}.
			\label{hel_cur8}
		\end{equation}
		Then, using equations (\ref{hel_cur5}), (\ref{hel_cur6}) and (\ref{hel_cur7}), we find
		\begin{equation}
			\ddot{x}=\frac{\varepsilon\vartheta}{\sqrt{|\varphi|}}\dot{V}+\frac{ \sqrt{|\varphi|}\kappa_\lambda'}{\vartheta
				\kappa_\lambda }\left(\frac{\vartheta^2}{\varphi}+\frac{\varepsilon a}{d}g(\dot{V},\dot{V})\right)V,
			\label{hel_cur9}
		\end{equation}
		hence 
		\begin{equation*}
			\nabla_{\dot{x}}\ddot{x}=-\frac{\varepsilon\vartheta}{\sqrt{|\varphi|}}g(\dot{V},\dot{V})V+\frac{\sqrt{|\varphi|}\kappa_\lambda'}{\vartheta
				\kappa_\lambda }\left(\frac{\vartheta^2}{\varphi}+\frac{\varepsilon a}{d}g(\dot{V},\dot{V})\right)\dot{V}.
		\end{equation*}
		Then, substituting from equations (\ref{hel_cur5}), (\ref{hel_cur7}), (\ref{hel_cur8}) into equations (\ref{sl_cur1}) and (\ref{sl_cur2}),
		using the fact that $\{V,\dot{V}\}$ is a frame along $x$ and taking the $V$-component of the first equation and the $\dot{V}$-component of the second one, we find 			
		\begin{equation}			
			\varepsilon\vartheta\left[-\varepsilon_\lambda\varepsilon_1\varepsilon_2\vartheta +\frac{1}{\kappa_\lambda\tau}\left(\varepsilon_\lambda\varepsilon_1\kappa_\lambda^2-\frac{d}{2\alpha}(\varepsilon_\lambda-\varepsilon\vartheta^2)\right)f_2\right] =-\varepsilon_1\varepsilon_2,
			\label{hel_cur10}
		\end{equation}
		\begin{equation}
			-\frac{(\kappa_\lambda')^2}{\kappa_\lambda^3\tau}\left(1+\frac{\varepsilon d\vartheta^2}{a\varphi g(\dot{V},\dot{V})}\right)f_2 -\varepsilon_\lambda\varepsilon_1\varepsilon_2\vartheta+\frac{\varepsilon_\lambda\varepsilon_1\kappa_\lambda}{\tau}f_2 +\frac{\varepsilon d(a+c-d)\vartheta^2}{2\kappa_\lambda\tau \alpha\varphi}f_2=0.
			\label{hel_cur12}
		\end{equation}
		Comparing the value of $\frac{\tau}{f_2}$ in (\ref{hel_cur10}) and (\ref{hel_cur12}) we get 
		\begin{equation}				
			\kappa_\lambda^2=	\varepsilon_1(\varepsilon_\lambda-\varepsilon\vartheta^2)\left[\frac{\varepsilon d^2\vartheta^2}{\alpha\varphi} +\frac{(\kappa_\lambda')^2}{\kappa_\lambda^2}\left(1+\frac{\varepsilon d\vartheta^2}{a\varphi g(\dot{V},\dot{V})}\right) \right].
			\label{hel_cur16}
		\end{equation}
		Since $\frac{\kappa_\lambda'}{\kappa_\lambda}$ is constant, then (\ref{hel_cur16}) entails $\kappa_\lambda$ is constant, i.e. 
		$\kappa_\lambda'=0$. Therefore, equations (\ref{hel_cur6}) and \eqref{hel_cur5} yield, respectively, $\sigma=\frac{\vartheta^2}{|\varphi|}$ and $
			\dot{x}=\frac{\varepsilon}{\sqrt{|\varphi|}}\vartheta V.$
			
		 On the other hand, equations (\ref{hel_cur16}) and \eqref{equat11} become, respectively,
		\begin{equation}\label{obl-kappa}
			\kappa_\lambda^2=\frac{\varepsilon_1d^2\vartheta^2}{\alpha\varphi}(\varepsilon\varepsilon_\lambda-\vartheta^2), \quad \textup{and} \quad f_2=\frac{\varepsilon_1d\vartheta}{2\alpha\varphi\kappa_\lambda\tau}(\varepsilon\varepsilon_\lambda-\vartheta^2)(\varepsilon\varphi-2\varepsilon_\lambda d\vartheta^2).
		\end{equation}  
		Substituting from the last two identities into (\ref{hel_cur12}), we find
		\begin{equation}\label{obl-tau}				
			\tau^2=\frac{\varepsilon_\lambda\varepsilon_1\varepsilon_2}{4\alpha\varphi}(\varepsilon\varphi-2\varepsilon_\lambda d\vartheta^2)^2.
		\end{equation}
		It follows from the inequality $\varepsilon_2(\varepsilon_\lambda-\varepsilon\vartheta^2)>0$ and \eqref{obl-kappa} that $\varepsilon_1\varepsilon_2 \alpha >0$. In a similar way, \eqref{obl-tau} yields $\varepsilon_\lambda\varepsilon_1\varepsilon_2\alpha\varphi >0$ which, together with $\varepsilon_1\varepsilon_2 \alpha >0$, gives $\varepsilon_\lambda\varphi>0$. Hence $\varepsilon_\lambda$ and $\varphi$ have the same sign and therefore $\varepsilon\varepsilon_\lambda=1$.
		
		Conversely, if an oblique Frenet curve $\lambda$ satisfies the five conditions of the theorem yield \eqref{cur_hel1} and \eqref{cur_hel2}, i.e. $\lambda$ is a helix curve with direction the geodesic vector field.
		\cqfd
		
		
		\subsection*{Cartan null helices in $T_1M^2(\kappa)$}
		
		
		Since we are concerned with null curves in $T_1M^2(\kappa)$, the manifold $(T_1M^2(\kappa),\tilde{G})$ must have a signature of either $(1,2)$ or $(2,1)$. Without loss of generality, we will assume that $T_1M^2(\kappa)$ is a Lorentzian manifold.

		We now provide a characterization of non-geodesic null helix curves in $T_1M^2(\kappa)$ with direction the geodesic vector field on  $T_1M^2(\kappa)$. Let $\lambda=(x,V)$ be a non-geodesic null curve in $T_1M^2(\kappa)$, then $\tilde{G}(\lambda'',\lambda'')> 0$. We parametrize $\lambda$ using the pseudo-arc parameter, that is such that $\tilde{G}(\lambda'',\lambda'')=1$. The Cartan Frenet frame of $\lambda$ satisfies the Cartan Frenet equations:
		\begin{equation}\label{cart-frenet}
			\left\{\begin{split}
				\tilde{\nabla}_TT=&W,\\
				\tilde{\nabla}_TW=&-\kappa_\lambda T-N,\\
				\tilde{\nabla}_TN=&\kappa_\lambda W,
			\end{split}\right.
		\end{equation}
		where  $T=\lambda'$, $W=\lambda''$, and $N$ is the transversal vector field given by $N=-\lambda^{(3)}-\frac{1}{2}\tilde{G}(\lambda^{(3)},\lambda^{(3)})\lambda'$ 
		and $\kappa_\lambda=\frac{1}{2}\tilde{G}(\lambda^{(3)},\lambda^{(3)})$ is the \emph{lightlike curvature} of $\lambda$. It is straightforward to show that: 
		\begin{equation}\label{null_hel01}
			\tilde{G}(T,T)=\tilde{G}(N,N)=\tilde{G}(T,W)=\tilde{G}(N,W)=0, \quad \tilde{G}(T,N)=1.
		\end{equation}
		Similar to the non-null case, let  $\vartheta=\tilde{G}(\lambda',\tilde{\xi})$. Then  $\vartheta=\frac{\varphi}{\sqrt{|\varphi|}}g(\dot{x},V),$
		hence  
		\begin{equation}
			g(\dot{x},V)=\frac{\varepsilon\vartheta}{\sqrt{|\varphi|}}.
			\label{null_hel}
		\end{equation}
		
		In the following theorem, we will demonstrate that any non-geodesic null helix curve in $T_1M^2(\kappa)$ directed along the geodesic vector field, satisfies a specific differential equation, which will be used to further characterize these curves.

		\begin{Theorem} \label{cartan_cur}
			Let $\lambda:I\subset\mathbb{R}\rightarrow T_1M^2(\kappa)$ be a non-geodesic null helix curve parameterized by the pseudo-arc parameter, with direction the geodesic vector field and denote by $\kappa_\lambda$ its lightlike curvature, $\varepsilon=\frac{\varphi}{|\varphi|}$ and $\vartheta:=\tilde{G}(\lambda^\prime,\tilde{\xi})$. Then $\lambda$ is a solution of the differential equation 
			\begin{equation}
				\sqrt{|\varphi|}\vartheta\lambda^{(3)}+\sqrt{|\varphi|}\left[\kappa_\lambda\vartheta-h_1\right]\lambda'+V^h=0,
				\label{null_sl_cur}
			\end{equation} 
			where 
			\begin{equation}
				h_1=\frac{\varepsilon\sqrt{|\varphi|}}{2}[g(\ddot{x},\dot{V})-g(\ddot{V},\dot{x})]-\vartheta g(\dot{V},\dot{V})-\frac{\varepsilon d}{2\alpha}\vartheta^3-\kappa_\lambda\vartheta,
				\label{null_sl_cur2}
			\end{equation}
			or equivalently, the base curve $x$ and the vector field $V$ along $x$ satisfy the system of differential equations
			\begin{equation}
				\begin{split} 
					\vartheta\nabla_{\dot{x}}\ddot{x}+\left[\kappa_\lambda\vartheta-h_1+\frac{d(a+c-d)\vartheta^3}{2\alpha|\varphi|}\right]\dot{x} +\frac{3\varepsilon(d-(a+c))}{2(a+c)\sqrt{|\varphi|}}\vartheta^2\ddot{V}+g(\dot{V},\dot{x})\vartheta\dot{V}&\\
					+\left[\frac{3\vartheta}{2}\dot{x}(g(\dot{V},\dot{x}))+\frac{\varepsilon \vartheta^2(3d-(a+c))}{2(a+c)\sqrt{|\varphi|}}g(\dot{V},\dot{V}) +\frac{d^2\vartheta^4}{\alpha\varphi\sqrt{|\varphi|}} +\frac{1}{\sqrt{|\varphi|}}
					\right] V&=0,
				\end{split}
				\label{lceq11}
			\end{equation}
			\begin{equation}
				\begin{split}
					\frac{\varepsilon(a+c-3d)\vartheta^2}{2a\sqrt{|\varphi|}}\ddot{x}+\vartheta\nabla_{\dot{x}}\ddot{V} +\frac{\varepsilon(a+c-d)\vartheta^2}{2a\sqrt{|\varphi|}}g(\dot{x},\dot{V})V&\\
					+\left[\frac{d\vartheta^3}{a|\varphi|}-\frac{(a+c-d)^2\vartheta^3}{2\alpha|\varphi|}-\vartheta g(\ddot{V},V)+\kappa_\lambda\vartheta-h_1\right]\dot{V}&=0.
				\end{split}
				\label{lceq2}
			\end{equation}
		\end{Theorem}
         \emph{Proof of Theorem \ref{cartan_cur}:}
			Since $\{ T,N,W\}$ is a local frame along $\lambda$, we can write $\tilde{\xi}$ as 
			$$\tilde{\xi}=h_1T+h_2 N+h_3W.$$
			The scalar product of $\tilde{\xi} $ by $T$ gives, by virtue of \eqref{null_hel01}, $h_2=\vartheta$ and since $\vartheta$ is constant, we have $$h_3=\tilde{G}(\tilde{\xi},W)=\tilde{G}(\tilde{\xi},\tilde{\nabla}_TT)=-\tilde{G}(\tilde{\nabla}_T\tilde{\xi},T) =-\varepsilon\tilde{G}(\tilde{\phi}(T),T)=-\varepsilon d\tilde{\eta}(T,T)=0.$$
			Then differentiating $\tilde{G}(\tilde{\xi},W)$, we obtain $\tilde{G}(\tilde{\nabla}_T\tilde{\xi},W)+\tilde{G}(\tilde{\xi},\tilde{\nabla}_TW)=0$ and, as a consequence, we get using again \eqref{null_hel01} and \eqref{helix-eq3}
			\begin{align*}
					h_1=&\tilde{G}(\tilde{\xi},N)=-\tilde{G}(\tilde{\xi},\lambda^{(3)})-\frac{1}{2}\tilde{G}(\lambda^{(3)},\lambda^{(3)})\vartheta\\
					=&\tilde{G}(\tilde{\nabla}_T\tilde{\xi},W)-\kappa_\lambda\vartheta\\
					=&\frac{\varepsilon\sqrt{|\varphi|}}{2}[g(\ddot{x},\dot{V})-g(\ddot{V},\dot{x})]-\vartheta g(\dot{V},\dot{V})-\frac{\varepsilon d}{2\alpha}\vartheta^3-\kappa_\lambda\vartheta.
			\end{align*}
			It follows then that $\tilde{\xi}=h_1T+\vartheta N$. Replacing $N$ by its expression, we get equation \eqref{null_sl_cur}.
			Now, we deduce from \eqref{vel-1} that 
			\begin{equation}\label{veloc-null}
				g(\dot{x},\dot{x})=-\frac{1}{a+c}\left[\frac{d\vartheta^2}{|\varphi|}+ag(\dot{V},\dot{V})\right].
			\end{equation}
			Then, decomposing equation (\ref{null_sl_cur}) into its horizontal and vertical parts and using \eqref{helix-eq5} and \eqref{veloc-null}, we find equations (\ref{lceq11}) and (\ref{lceq2}). 
		\cqfd
			\begin{Corollary} \label{legendre_cur}
		    If $\lambda$ a null Legendre curve in $T_1M^2(\kappa)$, then it is a vertical curve.
		\end{Corollary}
		\emph{Proof of Corollary \ref{legendre_cur}:}
			Suppose that $\lambda$ is a non-vertical Legendre curve ($\vartheta=0$). Then equations (\ref{lceq11}) and (\ref{lceq2}) are equivalents, respectively, to 
			\begin{equation*}
				\sqrt{|\varphi|}h_1\dot{x} -V=0 \quad \textup{and} \quad	h_1\dot{V}=0,
			\end{equation*}
			the secod equation implies $h_1=0$, then substituting into the first equation, we find $V=0$, which contradicts $V$ is being unitary.

		\emph{Proof of Theorem \ref{null-hor-hel}:}
		By assumption, we have $\dot{V}=0$, then the fact that  $\lambda$ is a lightlike curve implies 
		\begin{equation}
			g(\dot{x},\dot{x})=-\frac{d\vartheta^2}{(a+c)|\varphi|}.
			\label{null_hel1}
		\end{equation}
		Using the expression (\ref{Der2-Lambda}) of $\lambda''$ and the identity $\tilde{G}(\lambda'',\lambda'')=1$, we find 
		\begin{equation}
			g(\ddot{x},\ddot{x})=\frac{1}{a+c}\left[1+\frac{d^2\vartheta^4}{\alpha\varphi}\right].
			\label{null_hel12}
		\end{equation} 
		It follows that equations (\ref{lceq11}) and (\ref{lceq2}) become, respectively,
		\begin{equation}
			\vartheta\nabla_{\dot{x}}\ddot{x}+\left[\kappa_\lambda\vartheta-h_1+\frac{d(a+c-d)\vartheta^3}{2\alpha|\varphi|}\right]\dot{x}
			+\frac{1}{\sqrt{|\varphi|}}\left[\frac{d^2\vartheta^4}{\alpha\varphi} +1\right]V=0
			\label{null_hel2}
		\end{equation}
		and
		\begin{equation}
			\varepsilon(a+c-3d)\vartheta^2\ddot{x}=0,
			\label{null_hel3}
		\end{equation}
		
		where \begin{equation}
			h_1=-\frac{\varepsilon d}{2\alpha}\vartheta^3-\kappa_\lambda\vartheta.
			\label{null_hel4}
		\end{equation}	
		 By Corollary \ref{legendre_cur}, $\vartheta\neq 0$, then we deduce from equation (\ref{null_hel3}) that either  $a+c-3d=0$ or  $\ddot{x}=0$, that is $x$ is a geodesic.
		\begin{itemize}
			\item[\textbf{Case 1:}] If $\ddot{x}=0$, then equation (\ref{null_hel12}) implies that $d\neq 0$ and
			\begin{equation}\label{null_hel41}
				\vartheta^4=-\frac{\alpha\varphi}{d^2}.
			\end{equation} 
			Substituting from \eqref{null_hel41} into (\ref{null_hel2}) and using \eqref{null_hel4}, we get
			\begin{equation}\label{null_hel51}
				\kappa_\lambda = -\frac{\varepsilon d}{2a\varphi}\vartheta^2.
			\end{equation}
			Recalling that, according to the nature of the (para-)contact structure on $T_1M^2(\kappa)$, we have either $|\varphi|=4\alpha$ or $|\varphi|=-4\alpha$ (cf. Proposition \ref{contact}), we shall treat the two cases.
			\begin{itemize}
				\item If $|\varphi|=4\alpha$, then \eqref{null_hel41} yields $\varphi<0$, i.e. $\varepsilon <0$. Hence $\vartheta^2=-\frac{\varphi}{2|d|}$, then substituting into \eqref{null_hel51}, we find $\kappa_\lambda=\frac{ d(a+c)}{|d|\varphi}\vartheta^2$.
				\item If $|\varphi|=-4\alpha$, then \eqref{null_hel41} implies  $\varphi>0$, i.e. $\varepsilon >0$, by consequence $\vartheta^2=\frac{\varphi}{2|d|}$, then \eqref{null_hel51} gives $\kappa_\lambda=-\frac{d(a+c)}{|d|\varphi}\vartheta^2$.
			\end{itemize}
			\item[\textbf{Case 2:}] Let us suppose that $x$ is not a geodesic and $a+c=3d$, then $d\neq 0$ and hence $g(\dot{x},\dot{x})\neq 0$ and, from equation (\ref{null_hel12}), we get  $g(\ddot{x},\ddot{x})=\frac{1}{3d}\left[1+\frac{d\vartheta^4}{4\alpha}\right]$ is constant.
			\begin{itemize}
				\item If $\vartheta^4=-\frac{4\alpha}{d}$, that is  $\ddot{x}$ is a lightlike vector field along $x$, then the vector fields $\lbrace \dot{x},\ddot{x}\rbrace$ are linearly independent because of their causal character. We deduce that $\lbrace \dot{x},\ddot{x}\rbrace$ is a local frame along  $x$. 
				From $g(\ddot{x},\dot{x})=0$, we find 
					$g(\nabla_{\dot{x}}\ddot{x},\dot{x})=0,$
				thus 
				\begin{equation}
					\nabla_{\dot{x}}\ddot{x}=0.
					\label{null_hel9}
				\end{equation}
				In this case, we find the same result as in case 1, i.e. when $\ddot{x}=0$.
				\item If $\vartheta^4\neq-\frac{4\alpha}{d}$, then we put 
				 $\varepsilon'=\frac{g(\ddot{x},\ddot{x})}{|	g(\ddot{x},\ddot{x})|}$.  Using the same argument in \ref{proof}, we find 
				\begin{equation*}
					\nabla_{\dot{x}}\ddot{x}=\frac{|\varphi|}{d\vartheta^2}\left[1+\frac{d\vartheta^4}{4\alpha}\right]\dot{x},
				\end{equation*}
				that is $x$ is a pseudo-Riemannian circle in $M^2(\kappa)$. Therefore, equation (\ref{null_hel2}) becomes 
				\begin{equation*}
					A\dot{x}
					+BV=0,
				\end{equation*}
			where \begin{equation*}
					A=\frac{4\varepsilon}{\vartheta^2}\left[1+\frac{d\vartheta^4}{4\alpha}\right]+2\kappa_\lambda\vartheta+\frac{\varepsilon\vartheta^3}{4a}, \quad \textup{and}
				    \quad 
					B=\frac{1}{\sqrt{|\varphi|}}\left[\frac{ d\vartheta^4}{4\alpha}+1 \right].
				\end{equation*}
				Since $B \neq 0$, by hypothesis, then $A \neq 0$ and $
					V=-\frac{A}{B}\dot{x}.$
				Then we have $1=g(V,V)=\frac{A^2}{B^2}g(\dot{x},\dot{x})$, which yields, by virtue of \eqref{null_hel1},
				\begin{equation*}\label{null_hel93}
					\frac{A^2}{B^2}=-\frac{3|\varphi|}{\vartheta^2}< 0,
				\end{equation*}
				 which is a contradiction.
			\end{itemize}
		\end{itemize}
		
		The converse part of the theorem follows from the second condition of it.
		\cqfd
		

		\emph{Proof of Theorem \ref{null-obl-hel}:}
		Suppose that $g(\dot{x},\dot{x})$ is a constant $\sigma$. Since $g(\dot{x},V)=\frac{\varepsilon}{\sqrt{|\varphi|}}\vartheta$, then $g(\ddot{x},V)=-g(\dot{x},\dot{V})$,  taking the scalar product of (\ref{lceq2}) by $V$ we obtain
		\begin{equation}
			\begin{split}
			a\sqrt{|\varphi|}\vartheta g(\nabla_{\dot{x}}\ddot{V},V)+\varepsilon d\vartheta^2g(\dot{x},\dot{V})=0.
				\label{null_sl1}
			\end{split}
		\end{equation}
		 On the other hand, since $\lambda$ is a lightlike curve, then 	
			$$\sigma=-\frac{1}{a+c}\left(\frac{d\vartheta^2}{|\varphi|}+ag(\dot{V},\dot{V})\right),$$
		 which implies that $g(\dot{V},\dot{V})$ is constant 
		hence $g(\ddot{V},\dot{V})=0$. Differentiating twice the expression $g(\dot{V},V)=0$, we obtain $g(\nabla_{\dot{x}}\ddot{V},V)=-3g(\ddot{V},\dot{V})=0,$
		hence equation (\ref{null_sl1}) becomes
		\begin{equation*} d\vartheta^2g(\dot{x},\dot{V})=0.
		\end{equation*}	
		It follows that either $d=0$ or $\vartheta=0$ or $g(\dot{x},\dot{V})=0$. The scalar product (\ref{lceq11}) by $V$ implies that  $\vartheta \neq 0$.  
		Suppose $d=0$, then equation (\ref{lceq11})  becomes 
		\begin{equation}
			\vartheta\nabla_{\dot{x}}\ddot{x}+\left[\kappa_\lambda\vartheta-h_1\right]\dot{x}-\frac{3\varepsilon\vartheta^2}{2\sqrt{|\varphi|}}\ddot{V}+g(\dot{V},\dot{x})\vartheta\dot{V}
			-\frac{\varepsilon\vartheta^2}{2\sqrt{|\varphi|}}g(\dot{V},\dot{V}) V+\frac{1}{\sqrt{|\varphi|}}V=0.
			\label{null_sl4}
		\end{equation}
		Since $\{V,\dot{V}\}$ is an orthogonal frame along the curve $x$ and $g(\dot{V},\dot{V})\neq 0$, we have 
		\begin{equation}\label{tan_vector}
			\dot{x}=\frac{\varepsilon\vartheta}{\sqrt{|\varphi|}}V+\frac{g(\dot{x},\dot{V})}{g(\dot{V},\dot{V})}\dot{V},
		\end{equation}
		then the scalar product of \eqref{tan_vector} by $\dot{x}$ gives 
			 $$g(\dot{x},\dot{V})^2=\left(\sigma-\frac{\vartheta^2}{|\varphi|}\right)g(\dot{V},\dot{V}).$$ In particular $g(\dot{x},\dot{V})$ is a constant, that implies $g(\ddot{x},\dot{V})+g(\dot{x},\ddot{V})=0$, hence 
		\begin{equation}\label{tan_vector2}
			g(\ddot{x},\dot{V})-g(\dot{x},\ddot{V})=2g(\ddot{x},\dot{V}).
		\end{equation}
		Using $g(\dot{V},V)=0$, we find  $g(\ddot{V},V)=-g(\dot{V},\dot{V})$ and since $g(\ddot{V},\dot{V})=0$, then $
			\ddot{V}=-g(\dot{V},\dot{V})V,$ 
		hence equation \eqref{tan_vector} implies that 
		\begin{equation}\label{tan_vector3}
			\ddot{x}=-g(\dot{x},\dot{V})V+\frac{\varepsilon\vartheta}{\sqrt{|\varphi|}}\dot{V},
		\end{equation}
		thus $g(\ddot{x},\dot{V})=\frac{\varepsilon\vartheta}{\sqrt{|\varphi|}}g(\dot{V},\dot{V})$. Then substituting from \eqref{tan_vector2} into \eqref{null_sl_cur2}, we find that $h_1=-\kappa_\lambda\vartheta.$ 
		The scalar product of \eqref{null_sl4} by $\dot{V}$ gives then
		\begin{equation} \label{lig_obl_cur1}
			 g(\nabla_{\dot{x}}\ddot{x},\dot{V})+2\kappa_\lambda g(\dot{x},\dot{V})+ g(\dot{x},\dot{V})g(\dot{V},\dot{V})=0.
		\end{equation}
		Differentiating the expression \eqref{tan_vector3}, then taking the scalar product by $\dot{V}$, we obtain
	 $g(\nabla_{\dot{x}}\ddot{x},\dot{V})=-g(\dot{V},\dot{V})g(\dot{x},\dot{V})$. In this case, \eqref{lig_obl_cur1} 
		becomes \begin{equation*}
			 \kappa_\lambda g(\dot{x},\dot{V})=0,
		\end{equation*}
		thus either $\kappa_\lambda=0$ or $g(\dot{x},\dot{V})=0$.
		If $\kappa_\lambda=0$, since $\kappa_\lambda:=\frac{1}{2}\tilde{G}(\lambda^{(3)},\lambda^{(3)})$ then $\lambda^{(3)}$ is a lightlike vector field and $h_1=0$. In this case, equation (\ref{null_sl_cur}) becomes 
		\begin{equation}
			\sqrt{|\varphi|}\vartheta\lambda^{(3)}+V^h=0,
		\end{equation}
		which can not occur since $\tilde{G}(V^h,V^h)=\varphi\neq 0$, then $\kappa_\lambda\neq 0$. Hence $g(\dot{x},\dot{V})=0$, and consequently $\dot{x}=\frac{\varepsilon\vartheta}{\sqrt{|\varphi|}} V$, hence 
		\begin{equation} \label{tan_vector5}
			V=\frac{\varepsilon\sqrt{|\varphi|}}{\vartheta}\dot{x},\quad
			\sigma=\frac{\vartheta^2}{|\varphi|}\quad \textup{and}\quad g(\dot{V},\dot{V})=-\frac{\varepsilon}{a}\vartheta^2.
		\end{equation}
		It follows that  
		\begin{equation} \label{null_sl8}
			\dot{V}=\frac{\varepsilon\sqrt{|\varphi|}}{\vartheta}\ddot{x},\quad \textup{and} \quad
			\ddot{V}=\frac{\varepsilon\sqrt{|\varphi|}}{\vartheta}\nabla_{\dot{x}}\ddot{x}.
		\end{equation}
		Now, since $g(\dot{V},\dot{V})$ is constant, then we deduce from (\ref{null_sl8}) that $g(\ddot{x},\ddot{x})$ is constant equal to$
			g(\ddot{x},\ddot{x})=-\frac{\vartheta^4}{a\varphi}.$
		In particular, $x$ is not a geodesic and $g(\nabla_{\dot{x}}\ddot{x},\ddot{x})=0.$ We note that $\{\dot{x},\ddot{x}\}$ is an ortghogonal frame along $x$ and 
		\begin{equation*}
			g(\nabla_{\dot{x}}\ddot{x},\dot{x})=-g(\ddot{x},\ddot{x})=\frac{\vartheta^4}{a\varphi},
		\end{equation*}
		therefore
		\begin{equation}\label{tan_vector4}
			\nabla_{\dot{x}}\ddot{x}=\frac{\varepsilon}{a}\vartheta^2\dot{x}, \quad \ddot{V}=\frac{\sqrt{|\varphi|}}{a}\vartheta\dot{x} \quad \textup{and} \quad  \nabla_{\dot{x}}\ddot{V}=\frac{\sqrt{|\varphi|}}{a}\vartheta\ddot{x}.
		\end{equation}
		Using  equation (\ref{null_sl_cur2}), we find $
			h_1=-\frac{\varepsilon d}{2\alpha}\vartheta^3-\kappa_\lambda\vartheta.$
		After a straightforward calculation using \eqref{tan_vector5}, \eqref{null_sl8} and \eqref{tan_vector4}, we find that equation (\ref{lceq11}) and equation (\ref{lceq2})  are equivalents, respectively, to 
		\begin{equation}
			\kappa_\lambda=-\frac{\varepsilon}{2}\left(\frac{1}{\vartheta^2}+\frac{d\vartheta^2}{\alpha}\right), \quad \textup{and} \quad 	\kappa_\lambda=\frac{ d(d-3(a+c))}{4\alpha|\varphi|}\vartheta^2.
			\label{lig_curv1}
		\end{equation}
		On the other hand, using \eqref{Der2-Lambda}, we find that 
		\begin{equation*}
			\tilde{\nabla}_TT=\frac{d}{a+c}\ddot{x}^h+\frac{\varepsilon (a+c)}{a\sqrt{|\varphi|}}\vartheta\dot{x}^t,
		\end{equation*}
		hence $\tilde{G}(\tilde{\nabla}_TT,\tilde{\nabla}_TT)=-\frac{d^2\vartheta^4}{\alpha\varphi},$ but 
		 $\lambda$ is parameterized by the pseudo-arc parameter which entails
				$\vartheta^4=-\frac{\alpha\varphi}{d^2}.$
			In particular $\alpha\varphi<0$, 
			therefore $\varphi=-4\alpha$ and  $\vartheta^2=\frac{|\varphi|}{2|d|}$. As a consequence, the first equation of (\ref{lig_curv1})  
			becomes $\kappa_\lambda=-\frac{d}{4a|d|},$ then the second equation implies $d=a+c$ and consequently $\kappa_\lambda=\frac{|d|}{\varphi}$.
			
		The converse part of the theorem follows from the fourth condition of it.
		\cqfd


\begin{thebibliography}{1}
			\bibitem{Abb_Amr} M.T.K. Abbassi and N. Amri, On $g$-natural conformal vector fields on unit tangent bundles,
			\textit{Czech. Math. J.} \textbf{71} (2020) 1--35.
			
			\bibitem{Abb_18} M.T.K Abbassi, N. Amri, and G. Calvaruso, Kaluza--klein type Ricci solitons on unit tangent sphere
			bundles, \textit{Diff. Geom. Appl.} \textbf{59} (2018) 184--203.
			
			\bibitem{Abb-Bou-Cal} M.T.K Abbassi, K. Boulagouaz, and G. Calvaruso, On the Geometry of the Null Tangent Bundle of a Pseudo-Riemannian Manifold, \textit{Axioms} \textbf{12} (2023) 903.
			
			\bibitem{Abb_Cal} M.T.K. Abbassi and G. Calvaruso, The curvature tensor of $g$-natural metrics on unit tangent sphere bundles, \textit{Int. J. Contemp. Math. Sc.} \textbf{3} (2008) 245--258.
			
			\bibitem{Abb_Cal2} M.T.K. Abbassi and G. Calvaruso, $g$-natural contact metrics on unit tangent sphere
			bundles, \textit{Monatsh. Math.} \textbf{151} (2007) 89--109.
			
			\bibitem{Abb-Cal12} M.T.K. Abbassi and G. Calvaruso, $g$-natural metrics of constant curvature on unit tangent sphere bundles, \textit{Arch. Math. (Brno)} \textbf{48}(2) (2012) 81--95.
			
			\bibitem{Abb_Kow0} M.T.K Abbassi and O. Kowalski, On $g$-natural metrics with constant scalar curvature on unit tangent sphere bundles, in \textit{Topics in Almost Hermitian Geometry and related fields,  Proceedings of the International Conference in Honor of K. Sekigawa's 60th birthday}, (World Scientific 2005), pp. 1--29 .
			
			\bibitem{Abb_Kow} M.T.K Abbassi and O. Kowalski, On Einstein Riemannian g-natural metrics on unit tangent sphere bundles, \textit{Ann. Glob. Anal. Geom.} \textbf{38}(1) (2010), 11--20.

            \bibitem{Abb-Sar}K.M.T. Abbassi, M. Sarih, On some hereditary properties of Riemannian $g$-natural metrics on tangent bundles of Riemannian manifolds,\textit{ Diff. Geom. Appl.} \textbf{22} (2005) 19--47.
			
			\bibitem{Bar} M. Barros, General helices and a theorem of Lancret. \textit{Proc. Amer. Math. Soc.} \textbf{125}(5) (1997) 1503-1509.
			
			\bibitem{Loubeau} M. Benyounes, E. Loubeau, L. Todjihounde, Harmonic maps and Kaluza--Klein metrics on spheres, \textit{Rocky Mountain J. Math.} \textbf{42}(3) (2012) 791--821.
			
			\bibitem{Cal07} G. Calvaruso,  Einstein-like metrics on three-dimensional homogeneous Lorentzian manifolds, \textit{Geom. Dedicata} \textbf{127}(1) (2007) 99--119.
			
			\bibitem{Cal_15} G. Calvaruso and V. Marti\'{i}n-Molina, Paracontact metric structures on the unit tangent
			sphere bundle, \textit{Ann. Mat. Pura Appl.}, \textbf{194}(5) (2015) 1359--1380.
		
			
			\bibitem{Cal_Per} G. Calvaruso and D. Perrone, Contact pseudo-metric manifolds, \textit{ Diff. Geom. Appl.} \textbf{28}(5) (2010) 615--63.
			
			\bibitem{Sca_Her} A.J. Di Scala  and G. Ruiz-Hern\'{a}ndez, Higher codimensional Euclidean helix submanifolds, \textit{Kodai Math. J.}, \textbf{33}(2) (2010) 192-210.
			
			\bibitem{dug96} K.L. Duggal  and A. Bejancu,  \textit{Lightlike Submanifolds of Semi-Riemannian Manifolds and Applications} (Springer, 1996).
			
			\bibitem{Fer_Gim} A. Ferr\'{a}ndez, A. Gim\'{e}nez and L. Pascual, Null generalized helices in Lorentz--Minkowski spaces, \textit{J. Phys. A: Math. Gen.} \textbf{8243} (2002) 35--39.
			
			\bibitem{Fer_Gim2} A. Ferr\'{a}ndez, A. Gimenez and P. Lucas, Null helices in Lorentzian space forms, \textit{Int. J. Mod. Phys. A} \textbf{16}(30) (2001) 4845--4863.
			
			\bibitem{hay1} H.A. Hayden,  On a generalized helix in a Riemannian n-space, \textit{Proc. London Math. Soc.} \textbf{2}(1) (1931) 337--345.
			
			\bibitem{Hou_Sun} Z.H. HOU and L. SUN, Slant curves in the unit tangent bundles of surfaces, \textit{Int. Schol. Res. Not.} \textbf{2013} (2013).
			
			\bibitem{Ino_Lee} J.I. Inoguchi and S. Lee, Null curves in Minkowski 3-space, Int. Electr. J. Geom.  \textbf{1}(2) (2008) 40--83.
			
			\bibitem{Cho-Ino-Lee} J.T. Cho, J-I. Inoguchi and J-E. Lee, On slant curves in Sasakian 3-manifolds, Bull. Austr. Math. Soc. \textbf{74}(3) (2006) 359-367.
			
			\bibitem{Kar} H.B. Karadag and M. Karadag, Null generalized slant helices in Lorentzian space, \textit{Diff. Geom. Dyn. Syst.} \textbf{10} (2008) 178--185.
			
			\bibitem{Lee} J.E. Lee,  Slant curves and contact magnetic curves in Sasakian Lorentzian 3-manifolds, \textit{Symmetry}, \textbf{11}(6) (2019) 784.
			
			\bibitem{Nak_Zam} G. Nakova and S. Zamkovoy, Slant and Legendre null curves in 3-dimensional Sasaki-like almost contact $B$-metric manifolds, J. Geom. \textbf{112} (2021) 1--24.
			
			\bibitem{Spi70} M. Spivak, \textit{A comprehensive introduction to differential geometry. Vol. 4} (Publish or Perish, Incorporated, 1970).
			
			\bibitem{Ku_Zu} Z. K\"{u}\c{c}\"{u}karslan, M. Bekta\c{s} and H. \"{o}ztekin, Inclined Curves of Null Curves in and Their Characterizations, \textit{Int. J. Open Probl. Compt. Math} (2012) 5--2.
			
			\bibitem{Tak} T. Takahashi, Sasakian manifold with pseudo-Riemannian metrics \textit{T\^{o}hoku Math. J.} \textbf{21} (1969) 271--290.
			
			\bibitem{Tri} M.M. Tripathi, E. Kili\c{c}, S.Y. Perkta\c{s} and S. Kele\c{s}, Indefinite Almost Paracontact Metric Manifolds, \textit{Int. J. Math.  Math. Sci.} \textbf{2010}, Article ID 846195, 19 pages.
			
			\bibitem{Wood1} C.M. Wood, On the energy of a unit vector field, \textit{Geom. Dedicata} \textbf{64} (1997) 319--330.
		\end{thebibliography}
	\end{document}